\RequirePackage{fix-cm} 
\documentclass[smallextended]{svjour3}       
\smartqed  

\usepackage{calc,amsfonts,amssymb,amsmath,bm,url,color,theorem,graphicx,cite,shortcuts_OPT, bbm}
\usepackage{psfrag,subfigure,float,hyperref}
\usepackage{algorithm}
\usepackage{algorithmic}

\usepackage{fullpage}

\hypersetup{
  colorlinks   = true, 
  urlcolor     = blue, 
  linkcolor    = blue, 
  citecolor    = blue   
}

\usepackage{microtype}
\usepackage{algorithmic}

\usepackage{graphicx, amsfonts,mathabx}
\usepackage{booktabs} 
\usepackage{hyperref}
\usepackage{bm,shortcuts_OPT} 
\usepackage{colortbl}

\usepackage{pgfplots}
\pgfplotsset{compat=1.3}
\usepackage{environ}
\usepackage{tikz}
\usetikzlibrary{arrows,shapes,calc,tikzmark,backgrounds,matrix,decorations.markings}
\usepgfplotslibrary{fillbetween}

\usetikzlibrary{external}
\tikzsetexternalprefix{./tikz_gen_figures/}
\usetikzlibrary{shadows.blur}
\definecolor{asuorange}{rgb}{1,0.699,0.0625}
\definecolor{asured}{rgb}{0.598,0,0.199}
\definecolor{asuborder}{rgb}{0.953,0.484,0}
\definecolor{asugrey}{rgb}{0.309,0.332,0.340}
\definecolor{asublue}{rgb}{0,0.555,0.836}
\definecolor{asugold}{rgb}{1,0.777,0.008}

\definecolor{lavander}{cmyk}{0,0.48,0,0}
\definecolor{violet}{cmyk}{0.79,0.88,0,0}
\definecolor{burntorange}{cmyk}{0,0.52,1,0}

\providecommand{\env@tikzpicture@save@env}{}
\providecommand{\env@tikzpicture@process}{}

\newif\ifconfver
\confverfalse     

\newif\ifplainver
\plainvertrue

\newtheorem{Prop}{Proposition}
\newtheorem{Theorem}{Theorem}

\newtheorem{Lemma}{Lemma}
\newtheorem{Assumption}{Assumption}

\def\addlegendimage{\csname pgfplots@addlegendimage\endcsname}

\title{Accelerating Incremental Gradient Optimization with Curvature Information}

\author{Hoi-To Wai, Wei Shi, C\'{e}sar A. Uribe, Angelia Nedi\'{c}, Anna Scaglione\thanks{This work has been partially supported by the NSF grant
CCF-1717391 and CUHK Direct Grant \#4055113. H.-T.~Wai is with Dept.~of SEEM, Chinese University of Hong Kong, Hong Kong, W.~Shi was with Dept.~of EE, Princeton University, NJ, USA. A.~Nedi\'{c} and A.~Scaglione are with the School of ECEE, Arizona State University, Tempe, AZ, USA. C.~A.~Uribe is with LIDS, Massachusetts Institute of Technology, MA, USA. E-mails: \url{htwai@se.cuhk.edu.hk},\url{Angelia.Nedich@asu.edu}, \url{Anna.Scaglione@asu.edu}, \url{cauribe@mit.edu}.
Preliminary version of this work has been presented at the 55th Annual Allerton Conference on Communication, Control, and Computing in October, 2017 \cite{wai2017curvature}.}\\[.5cm]
\emph{In memory of Dr.~Wei Shi, a respected friend and talented scholar.}
}

\begin{document}
\maketitle

\begin{abstract}
This paper studies an acceleration technique for incremental aggregated gradient 
({\sf IAG}) method through the use of \emph{curvature} information for 
solving strongly convex finite sum optimization problems. 
These optimization problems of interest arise in large-scale learning applications. 
Our technique utilizes a curvature-aided gradient tracking
step to produce accurate gradient estimates incrementally using Hessian information. 
We propose and analyze two methods utilizing the new technique,
the curvature-aided IAG ({\sf CIAG}) method and the accelerated CIAG ({\sf A-CIAG}) method, 
which are analogous to gradient method and Nesterov's accelerated gradient method, respectively.
Setting $\kappa$ to be the condition number of the objective function, we prove the $R$ linear convergence rates of $1 - \frac{4c_0 \kappa}{(\kappa+1)^2}$ for the {\sf CIAG} method, and $1 - \sqrt{\frac{c_1}{2\kappa}}$ for the {\sf A-CIAG} method, where $c_0,c_1 \leq 1$ are constants inversely proportional to the distance between the initial point and the optimal solution. When the initial iterate is close to the optimal solution, the $R$ linear convergence rates match with the gradient and accelerated gradient method, albeit {\sf CIAG} and {\sf A-CIAG} operate in an incremental setting with strictly lower computation complexity.
Numerical experiments confirm our findings.
The source codes used for this paper can be found on \url{http://github.com/hoitowai/ciag/}.
\end{abstract}



\section{Introduction}
Consider a finite sum optimization problem
with $m$ component functions and a $d$-dimensional decision variable:
\beq \label{eq:opt}  \textstyle
\min_{ \prm \in \RR^d } ~ F(\prm) \eqdef \sum_{i=1}^m f_i ( \prm ) \eqs.
\eeq
Problem \eqref{eq:opt} is motivated by the empirical risk minimization model 
where we are learning a 
parameter $\prm$ from a finite set of data. 
The component function $f_i(\prm)$ 
represents the statistical mismatch between $\prm$ and the $i$th piece 
of data collected. 
The aim of problem \eqref{eq:opt} is to learn the
optimal parameter, denoted by $\prm^\star$, that fits with the $m$
available data points \cite{vapnik1999overview}. 

We are interested in the setting when each of the component function 
$f_i (\prm)$ is twice continuously differentiable and
the sum function $F(\prm)$ is strongly convex. Especially, this paper focuses on \emph{large-scale} instances of \eqref{eq:opt} with $m \gg 1$.
The difficulty with solving \eqref{eq:opt} lies in the overwhelming size of $m$,
which prohibits us from even applying simple 
first order methods.
For example, 
the \emph{full gradient} ({\sf FG}) 
requires the recursion: let $\gamma > 0$ be a step size,
\beq \label{eq:fg}  \textstyle
\prm^{k+1} = \prm^k - \gamma \sum_{i=1}^m \grd f_i ( \prm^k ) \eqs,
\eeq
that
involves a computation cost of ${\cal O}(md)$ 
floating point operations (FLOPS) per iteration to compute
the sum of $m$ gradient vectors.
As $m \gg 1$, this is   undesirable from a practical standpoint.  
To this end, a popular 
yet powerful approach is to adopt the so-called \emph{incremental} 
(or stochastic) methods
where only one of the component functions, e.g., the $i_k$th one, $f_{i_k}(\prm)$, 
is explored at the $k$th iteration. 
Examples include the incremental gradient ({\sf IG}) \cite{bertsekas2011incremental},
incremental aggregated gradient ({\sf IAG}) \cite{blatt2007convergent,vanli2016stronger,gurbuzbalaban2017convergence} methods when $i_k$ is deterministic; and 
the stochastic gradient ({\sf SG}) \cite{robbins1951stochastic}, stochastic average gradient 
({\sf SAG}) \cite{schmidt2017minimizing}, {\sf SAGA} \cite{defazio2014saga}, stochastic 
variance reduced gradient ({\sf SVRG}) \cite{xiao2014proximal} methods when
$i_k$ is chosen randomly. 
Related work can be found in \cite{mairal2015incremental,so2017non} 
for the case with non-convex functions; 
in \cite{lan2017optimal} for the cases of finite-sum and primal-dual optimization;
in \cite{nedic2001convergence,nedic2001incremental} for the case with non-smooth functions. The interested readers are referred to \cite[Section 4, 5]{bottou2016optimization} for a comprehensive overview on the topic.

A number of the above incremental methods achieve linear convergence 
via the variance reduction technique, e.g., \cite{xiao2014proximal,defazio2014saga,schmidt2017minimizing}.  
In the worst case, these methods require $k = {\cal O}( m \log (1/ \epsilon) )$ iterations to 
guarantee that they compute an $\epsilon$-optimal solution $\prm^k$ satisfying
$F(\prm^k) - F(\prm^\star) \leq \epsilon$. 
In problem instances when $m \gg 1$, such rates may not be ideal. In fact, recent work 
\cite{agarwal2015lower,arjevani2016dimension,lan2017optimal}
showed that the dependence on $m$
is necessary for incremental methods whose updates are 
linear combinations of the first order information.
To improve the convergence rate of incremental methods, 
a recent direction is to adopt ideas of second-order
optimization. Examples include \cite{gurbuzbalaban2015globally,rodomanov2016superlinearly,mokhtari2017iqn}
which extended Newton and 
quasi-Newton methods to the incremental setting, resulting in {\sf NIM} 
\cite{rodomanov2016superlinearly} and {\sf IQN} \cite{mokhtari2017iqn}. 
These works demonstrated that at the expense of additional storage
or computation cost,
it is possible to achieve local but \emph{superlinear} convergence.

We propose a \emph{curvature-aided gradient tracking} technique 
for accelerating incremental gradient methods using the curvature (Hessian) information. 
By applying Taylor expansion on the component functions' gradients, 
we derive a new gradient estimator whose error depends on the \emph{squares}
of the optimality gaps. Based on this new gradient tracking technique, 
we propose two incremental gradient methods, called
Curvature-aided Incremental Aggregated Gradient ({\sf CIAG}) method and 
Accelerated {\sf CIAG} ({\sf A-CIAG}) method.
The proposed methods require ${\cal O}(md)$ storage cost and ${\cal O}(d^2)$ 
computation cost per iteration. 
Furthermore, we characterize the $R$-linear convergence rates of the proposed methods. 
For the {\sf CIAG} (\resp {\sf A-CIAG}) method, we show that the rate is $1 - \frac{4 c_0 \kappa}{(\kappa+1)^2}$ (\resp $1 - \sqrt{\frac{c_1}{2 \kappa}}$), where $\kappa$ is the conditional number of the objective function of the summed objective function $F(\prm)$ in \eqref{eq:opt}, $c_0,c_1$ are constants which can be set to $1$ given that the initial iterate is close to the optimal solution [cf.~Theorems~\ref{thm:scvx} and \ref{thm:main}]. 
In addition, we establish the $R$-linear convergence rates
for two \emph{non-linear} inequalities [cf.~Propositions~\ref{prop:gen_ciag} and \ref{prop:gen_aciag}] which are of independent interest. 
In detail, our result reveals that the convergence of {\sf CIAG} and {\sf A-CIAG} methods depends on the initialization,
and the trajectory of convergence can be divided into two 
phases --- the initial phase that exhibits a slow but global linear rate;
and the asymptotic phase that exhibits linear convergence of an accelerated rate.  

A comparison of the proposed {\sf CIAG}/{\sf A-CIAG} methods to state-of-the-art methods is summarized in Table~\ref{tab:com}. Compared to first order methods such as {\sf SAG}, we remark that the proposed methods (as well as {\sf SVRG2}, {\sf IQN}, {\sf NIM}) require the Hessian of the objective function to be Lipschitz continuous. Though the storage requirement for the proposed methods is ${\cal O}(md+d^2)$ which is higher than the ${\cal O}(md)$ requirement of the {\sf SAG} methods when $d \gg 1$, the number of iterations required is significantly lower when the desired accuracy is high, $\epsilon \ll 1$. Overall, as seen from Table~\ref{tab:com}, the proposed methods are desired when $d$ is small-to-moderate, and the number of component functions $m$ is large.

\begin{table}[t]
\begin{center}
\renewcommand\baselinestretch{1.2}\selectfont
{ 
\begin{tabular}{l l l l} 
\toprule 
& Storage & Comp. & 
\# iterations to $\epsilon$-optimal solution
\\
\midrule 
{\sf IAG} \cite{blatt2007convergent} & ${\cal O}(md)$ & ${\cal O}(d)$ & ${\cal O}( m\!~ \kappa \log(1/\epsilon) )$ {\small [worst-case]} \\
\hline
{\sf SAG} \cite{schmidt2017minimizing} & ${\cal O}(md)$ & ${\cal O}(d)$ & ${\cal O}(\max\{ \kappa, m \} \log(1/\epsilon) )$ {\small[expected]}\\
\hline
{\sf SAGA} \cite{defazio2014saga} & ${\cal O}(md)$ & ${\cal O}(d)$ & ${\cal O}( ( \kappa + m) \log(1/\epsilon) )$ {\small[expected]}\\
\hline
{\sf SVRG} \cite{xiao2014proximal} & ${\cal O}(d)$ & ${\cal O}(d)$ & ${\cal O}( ( \kappa + m) \log(1/\epsilon) )$ {\small[expected]} \\
\hline
{\sf AccSVRG} \cite{nitanda2014stochastic} & ${\cal O}(d)$ & ${\cal O}(d)$ & ${\cal O}( ( 
m \sqrt{ \kappa } + m
 ) \log(1/\epsilon) )$ {\small[expected]} \\
\hline
\hline
{\sf SVRG2} \cite{gower2017tracking} & ${\cal O}( d^2 )$ & ${\cal O}(d^2)$  & ${\cal O}( \kappa \log(1/\epsilon) + m )$  {\small[expected]}\\
\hline
\hline
{\sf IQN} \cite{mokhtari2017iqn} & ${\cal O}(md^2)$ & ${\cal O}(d^2)$ & ${\cal O}(m)$,~\ie super-linear {\small[worst-case]}\\ 
\hline
{\sf NIM} \cite{rodomanov2016superlinearly} & ${\cal O}(md+d^2)$ & ${\cal O}(d^3)$ & ${\cal O}(m)$,~\ie super-linear {\small[worst-case]}\\
\hline
\hline
\cellcolor{black!10}{\sf CIAG} (proposed) & \cellcolor{black!10}${\cal O}(md+d^2)$ & \cellcolor{black!10}${\cal O}(d^2)$ & \cellcolor{black!10}${\cal O}({\kappa} \log(1/\epsilon) + m )$ {\small[worst-case]}\\
\hline
\cellcolor{black!10}{\sf A-CIAG} (proposed) & \cellcolor{black!10}${\cal O}(md+d^2)$ & \cellcolor{black!10}${\cal O}(d^2)$ & \cellcolor{black!10}${{\cal O}(\sqrt{\kappa} \log(1/\epsilon) + m )}$ {\small[worst-case]} \\
\bottomrule
\end{tabular}}\vspace{.2cm}
\end{center}
\caption{Comparison of state-of-the-art methods
for the strongly convex problem \eqref{eq:opt}. 
The second column is the memory required for the working variables.  
The third column is the \emph{per iteration} computation complexity in FLOPS. 
The last column is the (expected or worst-case) number of iterations to reach an $\epsilon$-optimal solution.
{\sf SVRG}, {\sf AccSVRG}, {\sf SVRG2}  
require recomputing the full gradient at every epoch and a careful tuning of the epoch size. 
The constant $\kappa$ is the condition number of $F(\prm)$ [see \eqref{eq:smooth_l} and \eqref{eq:strcvx} 
for the definition].
The rates for the last four methods are \emph{asymptotic}, \ie
they hold only when the desired accuracy is small, e.g., $\epsilon \ll 1$.
The shaded rows denote the proposed methods which achieve a good trade-off between per-iteration complexity and convergence speed.} \label{tab:com} \vspace{-.4cm}
\end{table}

It is worth mentioning that recently in \cite{gower2017tracking}
the authors proposed a method that incorporates
Hessian information to accelerate 
{\sf SVRG} method, giving the {\sf SVRG2} method. 
The authors developed approximation techniques to reduce the per iteration complexity 
from ${\cal O}(d^2)$ to ${\cal O}(d)$. 
The best convergence rates achieved therein
match that of the {\sf FG} method and is the same as the {\sf CIAG} method.
However, we note that {\sf SVRG2} re-computes the full Hessian
at the beginning of each epoch. This costs an additional ${\cal O}(m d^2)$ FLOPS per epoch
and the cost  is negligible only when the epochs are long, e.g., when the 
epoch lengths are ${\cal O}(m)$. 

 
\vspace{.2cm}
\noindent \textbf{Organization.}  Section~\ref{sec:back} studies 
incremental algorithms with \emph{gradient tracking} and introduces the curvature-aided gradient tracking technique. Section~\ref{sec:alg}
describes the proposed {\sf CIAG} and {\sf A-CIAG} methods, and discusses the  implementation issues.
Section~\ref{sec:ana} provides the main convergence results. 
Section~\ref{sec:num} demonstrates the efficacy of 
the proposed methods via numerical experiments. 

\vspace{.2cm}
\noindent \textbf{Notations.}
We use ${\bm 0}$ to denote an all-zero vector/matrix with suitable dimensions. 
Unless otherwise specified, we denote by $\| \cdot \|$ the Euclidean norm. 
A function $f : \RR^d \rightarrow \RR$ is $L$-smooth if 
\beq \label{eq:smooth_l}
f( \prm' ) \leq f( \prm ) + \langle \grd f( \prm ), \prm' - \prm \rangle + (L/2) \| \prm' - \prm \|^2 \eqs,
\eeq
for $\prm,\prm' \in \RR^d$,
and 
it is $\mu$-strongly convex if 
\beq \label{eq:strcvx}
f( \prm' ) \geq f( \prm ) + \langle \grd f( \prm ), \prm' - \prm \rangle + (\mu/2) \| \prm' - \prm \|^2 \eqs,
\eeq
for $\prm,\prm' \in \RR^d$.
Define $\kappa(f) \eqdef L / \mu$ as the condition number of $f$. 
Also, $f$ has an $L_H$-Lipschitz continuous Hessian if
\beq
\| \grd^2 f ( \prm' ) - \grd^2 f ( \prm ) \| \leq L_H \| \prm' - \prm \|,~\forall~\prm, \prm' \in \RR^d \eqs,
\eeq
where the norm on the left hand side is the matrix
norm induced by Euclidean norm. 
For a non-negative scalar sequence $\{ a^{(k)} \}_{k \geq 1}$, we say that 
it converges $R$-linearly at a rate $\rho$ if the sequence satisfies 
$\lim_{k \rightarrow \infty} a^{(k+1)} / a^{(k)} = \rho$, where $0 \leq \rho < 1$.
We use standard  Bachmann-Landau notations for 
asymptotic quantities: $a^{(k)} = {\cal O} ( f^{(k)} )$ (\resp 
$a^{(k)} = \Omega( g^{(k)} )$) implies that there exists 
$k_0 \in \ZZ_+$ and 
non-negative constant $C_1$ (\resp $C_2$) such that $a^{(k)} \leq C_1 f^{(k)}$
(\resp $a^{(k)} \geq C_2 g^{(k)}$) for all $k \geq k_0$. 
 
\section{Gradient Tracking in Incremental Methods} \label{sec:back}
We provide a high-level description of how to incorporate curvature information to accelerate incremental aggregated gradient method.
To set up the notation, let $i_k \in \{1,...,m\}$ be the component function index 
selected at the $k$th iteration, e.g., a simple rule is to use the cyclic rule as
 $i_k = (k~{\rm mod}~m)+1$, or  we may choose 
$i_k \sim {\cal U}\{1,...,m\}$. Let us define: 
\beq \label{eq:delay}
\tau_i^k \eqdef \max \{ \tau : \tau \leq k,~ i_\tau = i \} \eqs,
\eeq
\ie $\tau_i^k$ is the iteration index where the $i$th component function is last 
accessed after the completion of the $k$th iteration. 
As we focus  on analyzing the worst-case performance, we assume
that $\tau_i^k \in [k-K+1,k]$ for a constant $K = {\cal O}(m)$. 
For example, $K=m$ when $i_k$ is chosen 
by the cyclic rule.  

\vspace{.2cm}
\noindent \textbf{First-order Approximation.}
As described in \eqref{eq:fg}, at the $k$th iteration 
the gradient method computes the \emph{full gradient} 
as the complete sum $\sum_{i=1}^m \grd f_i( \prm^k )$. Such vector is unavailable 
in the incremental setting as only the access to the $i_k$th function is 
desired. A simple idea is to replace the full gradient by:\vspace{-.1cm}
\beq \label{eq:iag} \textstyle
  {\bm g}_{\sf IAG}^k \eqdef \sum_{i=1}^m \grd f_i ( \prm^{\tau_i^k} ) \approx \sum_{i=1}^m \grd f_i ( \prm^k ) \eqs,
\eeq
or equivalently using the recursion:\vspace{-.1cm}
\beq \label{eq:iag_recur} 
{\bm g}_{\sf IAG}^k = {\bm g}_{\sf IAG}^{k-1} - \grd f_{i_k} ( \prm^{\tau_{i_k}^{k-1}} ) + \grd f_{i_k} (\prm^k)  \eqs,
\eeq 
with ${\bm g}_{\sf IAG}^0 = \sum_{i=1}^m \grd f_i(\prm^0)$. 
The expression \eqref{eq:iag} describes the 
{\sf SAG} 
and {\sf IAG} methods, where their only differences lie
in the selection strategy for $i_k$. 

When the functions $f_i$ are smooth, we have
$\| \grd f_i ( \prm^{\tau_i^k} ) - \grd f_i ( \prm^k )  \| = {\cal O} ( \| \prm^{\tau_i^k} - \prm^k \| )$. 
Furthermore, 
\cite{gurbuzbalaban2017convergence} shows that the error is bounded:
\beq \label{eq:iag_err} {\textstyle
\| 
{\bm g}_{\sf IAG}^k - \sum_{i=1}^m \grd f_i (\prm^k)
 \|} = {\cal O} \big( \gamma m \max_{i=1,...,m} \| \prm^\star - \prm^{\tau_i^k} \| \big)  .
\eeq
Note that this error represents the `variance' in gradient estimation.
Eq.~\eqref{eq:iag_err} implies that the tracking error decays to zero
as long as $\prm^k$ converges to $\prm^\star$ [recall that $k - \tau_i^k \leq K < \infty$
and $K = \Theta(m)$]. 
However, 
the dependency on $m$ of the 
right hand side in \eqref{eq:iag_err} is undesirable as it leads to the
following   estimates:
\beq \label{eq:dyn_iag} \begin{split}
\| \prm^{k+1} - \prm^\star \|^2 & \leq \big(1- {\cal O}( 
\text{\scriptsize$\frac{\gamma}{\kappa(F)}$} ) \big) \| \prm^{k} - \prm^\star \|^2 + {\cal O}( \gamma^2 m^2 \hspace{-.2cm} \max_{ (k-2K)_{++} \leq \ell \leq k } \hspace{-.2cm} \| \prm^\ell - \prm^\star \|^2 ) \eqs,
\end{split}
\eeq
where $(x)_{++} \eqdef \max\{1, x\}$.
As analyzed by \cite{gurbuzbalaban2017convergence, feyzmahdavian2014delayed}, due to the dependence of $m^2$ on the right hand side,
the sequence of squared norm $\{ \| \prm^{k} - \prm^\star \|^2 \}_{k \geq 1}$ 
converges only when $\gamma = {\cal O}(1/m)$,
and finally this shows that $\| \prm^k - \prm^\star\|^2$
converges linearly at rate  
$1 - {\cal O}(1/ (m^2 \kappa(F)))$.   
Notice that a recent work \cite{vanli2016stronger} strengthened 
this rate to $1 - {\cal O}(1 / (m \kappa(F)))$. When the component function $i_k$ is selected independently at random for each $k$, the rate may also be improved,  see \cite{schmidt2017minimizing}. 

\vspace{.2cm}
\noindent \textbf{Second-order Approximation.}
The undesirable $m$-dependence of convergence rate with \eqref{eq:iag} is partly due to the crude approximation $\grd f_i( \prm^k ) \approx \grd f_i( \prm^{\tau_i^k})$ used. As a natural idea to improve the approximation, we consider the Taylor expansion applied on the $i$th \emph{gradient} vector itself:
\beq \label{eq:taylor}
\grd f_i ( \prm^k ) \approx \grd f_i ( \prm^{\tau_i^k} ) +
\grd^2 f_i ( \prm^{\tau_i^k} ) ( \prm^k - \prm^{\tau_i^k} ) \eqs.
\eeq
Compared to a first order approximation with $\grd f_i( \prm^k ) \approx \grd f_i( \prm^{\tau_i^k})$, the approximation given by the right hand side in \eqref{eq:taylor} has two important features. 
First, we observe that this term yields a \emph{second-order approximation} to the gradient of the component function $f_i(\prm)$ evaluated at $\prm^k$, and thus leading to an improved approximation.
Second, the approximate term on the right hand side does not require access to the $i$th function evaluated at $\prm^k$, which allows one to derive a similar recursion rule as \eqref{eq:iag_recur} to \emph{incrementally} compute the approximation. 
 
To further evaluate the approximation quality given by \eqref{eq:taylor}, for each $i$ we shall assume that the Hessian of $f_i$ is Lipschitz continuous.
We observe the following upper bound on the error:
\begin{Lemma} \label{lem:nes} \cite[Lemma 1.2.5]{nesterov2013introductory}
Assume that $\grd^2 f_i (\prm)$ is $L_{H,i}$ Lipschitz. Then:
\beq
\begin{split}
& \big\| \grd f_i( \prm ) - \big( \grd f_i ( \prm' ) + \grd^2 f_i ( \prm') ( \prm - \prm' ) \big) \big\|  \leq (L_{H,i}/2) \| \prm - \prm' \|^2,~\forall~\prm, \prm' \in \RR^d \eqs.
\end{split}
\eeq
\end{Lemma}
From Lemma~\ref{lem:nes}, 
we notice that the approximation error of the right hand side in  \eqref{eq:taylor}
depends on the \emph{squared difference} 
between $\prm^k$ and $\prm^{\tau_i^k}$. 
When $\prm^k$ is close to $\prm^{\tau_i^k}$, this error will be significantly
smaller than what is obtained for 
the first order approximation methods in \eqref{eq:iag}.
We remark that similar approximation schemes are proposed in recent work  \cite{zheng2017asynchronous,gower2017tracking}, which demonstrated the benefits of applying higher-order approximation in stochastic optimization. 

We refer to this form of gradient approximation in \eqref{eq:taylor} as the \emph{curvature-aided gradient tracking technique}. In the remainder of this paper, we shall develop and analyze practical optimization methods using the technique.
 \vspace{-.2cm}
  
\begin{algorithm}[t]
\caption{{\sf CIAG} and {\sf A-CIAG} Method.}\label{alg:aciag}
  \begin{algorithmic}[1]
  \STATE \textbf{Input}: Initial point $\prm^1 \in \RR^d$, step size parameter $\gamma > 0$, (for {\sf A-CIAG} only) extrapolation constant $\alpha$.
  \STATE Initialize the vectors/matrices in the memory:\vspace{-.1cm}
  \beq \label{eq:ciag_n_init}
   \prm_i^1 = \prm^1,~\tau_i^1 = 0,~\forall~i,~ {\bm b}^0 = {\bm 0},~{\bm H}^0 = {\bm 0} \eqs. \vspace{-.5cm}
  \eeq
  \FOR {$k=1,2,\dots$}
  \STATE \label{ciag_n:sel} Select $i_k \in \{1,...,m\}$, e.g., $i_k = (k~{\rm mod}~m) + 1$,
  and the algorithm is given access to $f_{i_k}(\prm)$.
  Set the counter variables as $\tau_{i_k}^k = k$, $\tau_j^k = \tau_j^{k-1}$ for all $j \neq i_k$.
  \STATE  \emph{(For {\sf A-CIAG} only)} If $k = 1$, we take $\eprm^{1} = \prm^1$; otherwise we compute the extrapolated variable as
 $ \eprm^k = \prm^k + \alpha ( \prm^k - \prm^{k-1} )$. 
  \STATE Update the vector in the memory:\vspace{-.2cm}
  \beq
  \prm_{i_k}^k = \begin{cases}
\prm^k, & \text{for {\sf CIAG} \eqs,} \\
\eprm^k, & \text{for {\sf A-CIAG} \eqs,} 
\end{cases}
  \eeq 
  and set $\prm_{j}^k = \prm_j^{k-1}$ for all $j \neq i_k$. 
   \STATE \label{ciag_n:sums} 
   If $\tau_{i_k}^{k-1} = 0$, then update: (\emph{self-initialization})\vspace{-.1cm}
   \beq \label{eq:ciag_n_inc_init}
   \begin{split}
   {\bm b}^{k} & = {\bm b}^{k-1} + \grd f_{i_k} (\prm_{i_k}^k) - \grd^2 f_{i_k} (\prm_{i_k}^k) \prm_{i_k}^k \eqs, \\
   {\bm H}^{k} & = {\bm H}^{k-1} + \grd^2 f_{i_k} (\prm_{i_k}^k) \eqs.\\[-.3cm]
   \end{split}
   \eeq
   Otherwise, update:\vspace{-.1cm}
   \beq \label{eq:ciag_n_inc}
   \begin{split}
   {\bm b}^{k} & = {\bm b}^{k-1} - \grd f_{i_k}( \prm_{i_k}^{k-1} ) + \grd f_{i_k} (\prm_{i_k}^k)  + \grd^2 f_{i_k}( \prm_{i_k}^{k-1} ) \prm_{i_k}^{k-1} - \grd^2 f_{i_k} (\prm_{i_k}^k) \prm_{i_k}^k \eqs, \\
   {\bm H}^{k} & = {\bm H}^{k-1} - \grd^2 f_{i_k}( \prm_{i_k}^{k-1} ) + \grd^2 f_{i_k} (\prm_{i_k}^k) \eqs.\\[-.3cm]
   \end{split}
   \eeq
   \STATE \label{ciag_n:upd} Compute the update:\vspace{-.1cm}
   \beq \label{eq:ciag_n_imp} \textstyle
   \prm^{k+1} = \prm_{i_k}^k - \gamma \big( {\bm b}^k + {\bm H}^k \prm_{i_k}^k   \big) \eqs. \vspace{-.1cm}
   \eeq
\ENDFOR
\STATE \textbf{Return}: an approximate solution to \eqref{eq:opt}, $\prm^{k+1}$.
  \end{algorithmic}
\end{algorithm}

\section{Proposed {\sf CIAG} and {\sf A-CIAG} Methods}
We propose two incremental methods, {\sf CIAG} and {\sf A-CIAG} methods, 
which are built using the curvature-aided gradient tracking technique \eqref{eq:taylor}. 
\label{sec:alg}
Our first method, the {\sf CIAG} method, performs the recursion:
\beq
\prm^{k+1} = \prm^k - \gamma {\bm g}_{\sf CIAG}^k \eqs,
\eeq
where $\gamma > 0$ is a step size and the gradient surrogate is given by
\beq \label{eq:ciag}
\begin{split}
{\bm g}_{\sf CIAG}^k &  \eqdef \sum_{i=1}^m \big( \grd f_i( \prm^{\tau_i^k} ) - \grd^2 f_i( \prm^{\tau_i^k} ) \prm^{\tau_i^k} \big) + \Big( \sum_{i=1}^m \grd^2 f_i( \prm^{\tau_i^k} ) \Big) \prm^k  \eqs,
\end{split}
\eeq
where $\tau_i^k$ was defined in \eqref{eq:delay}.
The method can be interpreted as applying the curvature-aided
gradient tracking [cf.~\eqref{eq:taylor}] to \emph{each} individual component function
and applying the {\sf FG} update \eqref{eq:fg}.  
On the other hand, the accelerated {\sf CIAG} ({\sf A-CIAG}) method
follows the recursion:
\begin{subequations} \label{eq:afg}
\begin{align} 
\eprm^k & = \prm^k + \alpha( \prm^k - \prm^{k-1} ) \eqs, \label{eq:afg_m} \\
\label{eq:afg_u} \prm^{k+1} &   = \eprm^k - \gamma {\bm g}_{\sf ACIAG}^k \eqs,
\end{align}
\end{subequations}
where $\alpha \in [0,1), \gamma > 0$ are predefined parameters
and the gradient surrogate is
\beq \label{eq:aciag} 
\begin{split}
{\bm g}_{\sf ACIAG}^k & \textstyle \eqdef \sum_{i=1}^m \big( \grd f_i( \eprm^{\tau_i^k}) - \grd^2 f_i( \eprm^{\tau_i^k})\eprm^{\tau_i^k} \big) + \Big( \sum_{i=1}^m \grd^2 f_i ( \eprm^{\tau_i^k} ) \Big) \eprm^k \eqs, \\[-.0cm]
\end{split}
\eeq
where $\tau_i^k$ was defined in \eqref{eq:delay}. 
The term $\eprm^k$ in \eqref{eq:afg_m} is called the extrapolated iterate,
which incorporates the `inertia' from previous iterates into the updates. 
Similar to the {\sf CIAG} method, 
the {\sf A-CIAG} method applies curvature-aided gradient tracking
to the gradient of each individual component function, evaluated at $\eprm^k$. 
 We remark that 
even though Hessians are used in the {\sf CIAG} and {\sf A-CIAG} methods, 
we do not attempt to compute any matrix inverses which is in contrast with Newton methods.
The Hessians here are  used to accelerate the gradient
tracking. 

The {\sf CIAG} and {\sf A-CIAG} methods are similar since
the full gradients evaluated at $\prm^k$ and $\eprm^k$ are approximated by their respective
curvature-aided gradient tracking approximates [cf.~\eqref{eq:ciag} \& \eqref{eq:aciag}]. 
Both methods can be implemented in a memory efficient
 incremental fashion. To see this, define 
\beq
\prm_i^k \eqdef \begin{cases}
\prm^{\tau_i^k}, & \text{for {\sf CIAG} \eqs,} \\
\eprm^{\tau_i^k}, & \text{for {\sf A-CIAG} \eqs,}
\end{cases}~~
{\bm H}^k \eqdef \sum_{i=1}^m \grd^2 f_i( \prm_i^k ),
\eeq  
and 
\beq \label{eq:bkhk} \textstyle
{\bm b}^k \eqdef \sum_{i=1}^m \big( \grd f_i( \prm_i^k ) - \grd^2 f_i( \prm_i^k ) \prm_i^k \big) \eqs.
\eeq
Note that ${\bm b}^k$ and ${\bm H}^k$ are the accumulated staled
iterates, gradients and Hessians. These can be computed incrementally through 
storing the staled iterates in memory.
Moreover, we can compute the gradient surrogates 
as ${\bm g}_{\sf CIAG}^k = {\bm b}^k + {\bm H}^k \prm^k$
and ${\bm g}_{\sf ACIAG}^k = {\bm b}^k + {\bm H}^k \eprm^k$. 
A pseudo code for implementing {\sf CIAG} and {\sf A-CIAG}  is provided
in Algorithm~\ref{alg:aciag}.\vspace{-.2cm}

\subsection{Computation and Storage Costs} 
We comment on the computation and storage costs of {\sf CIAG} and {\sf A-CIAG}. 
Note that \eqref{eq:ciag_n_inc_init}, \eqref{eq:ciag_n_inc} and \eqref{eq:ciag_n_imp} in the algorithm 
require ${\cal O}(d^2)$ FLOPS
and they are the dominant computation 
steps involved.
The overall complexities for {\sf CIAG} and {\sf A-CIAG} are 
thus ${\cal O}(d^2)$ FLOPS per iteration. 
On the other hand, the algorithm requires storing $m$ vectors of 
\mbox{$d$-dimesion} and a $d \times d$ matrix [cf.~\eqref{eq:ciag_n_init}], therefore the 
storage cost is ${\cal O}(md + d^2) = {\cal O}(md)$ if $m \geq d$. 
When $d$ is small, the computation and storage cost of {\sf CIAG} and {\sf A-CIAG} are
comparable to existing methods such as {\sf SAG}, {\sf SVRG}; 
meanwhile  for  
large $d$, the computation cost will be undesirable
for {\sf CIAG} and {\sf A-CIAG}.

When the component functions are the negative log-likelihood
of a linear model, the storage complexity can be reduced to 
as low as ${\cal O}(m)$. Note
linear models are common in machine learning problems.
We write\vspace{-.1cm}
\beq \label{eq:linearmodel}
f_i( \prm ) = g_i ( \langle \prm, {\bm x}_i \rangle ) + (\rho/2) \| \prm \|^2 \eqs,\vspace{-.1cm}
\eeq
where ${\bm x}_i$ represents the $i$th associated data, while $g_i : \RR \rightarrow \RR$
is twice continuously differentiable. 
Observe that\vspace{-.1cm}
\beq
\begin{split}
& \grd f_i ( \prm ) = g_i' ( \langle \prm, {\bm x}_i \rangle  ) {\bm x}_i + \rho \!~ \prm,~ \grd^2 f_i ( \prm ) = g_i'' ( \langle \prm, {\bm x}_i \rangle  ) {\bm x}_i ({\bm x}_i)^\top + \rho {\bm I} \eqs,
\end{split}\vspace{-.1cm}
\eeq
where ${\bm I}$ is the $d \times d$ identity matrix.
Substituting the above into \eqref{eq:ciag_n_inc} gives:
   \beq  
   \begin{split}
   {\bm b}^{k} & = {\bm b}^{k-1} + ( g_{i_k}' ( \langle \prm_{i_k}^k, {\bm x}_{i_k} \rangle)
   - g_{i_k}' (\langle \prm_{i_k}^{k-1}, {\bm x}_{i_k} \rangle) ) {\bm x}_{i_k} \\
   & +  g_{i_k}'' ( \langle \prm_{i_k}^{k-1}, {\bm x}_{i_k} \rangle) \langle \prm_{i_k}^{k-1}, {\bm x}_{i_k} \rangle {\bm x}_{i_k}  - g_{i_k}'' (\langle \prm_{i_k}^k, {\bm x}_{i_k} \rangle) \langle \prm_{i_k}^k, {\bm x}_{i_k} \rangle {\bm x}_{i_k} \vspace{.1cm}\\
   \end{split}
   \eeq
   \beq
      {\bm H}^{k} = {\bm H}^{k-1} +
   ( g_{i_k}'' ( \langle \prm_{i_k}^k, {\bm x}_{i_k} \rangle)
   - g_{i_k}'' (\langle \prm_{i_k}^{k-1}, {\bm x}_{i_k} \rangle) ) {\bm x}_{i_k} {\bm x}_{i_k}^\top \eqs.
   \eeq
Therefore, it suffices for {\sf CIAG} and {\sf A-CIAG} 
to keep
$\{ \langle {\bm x}_i, \prm_{i} \rangle \}_{i=1}^m$ for implementing the incremental updates, leading to an ${\cal O}(d^2+m)$ 
storage cost. 

\section{Convergence Analysis} \label{sec:ana}
In this section, we demonstrate that the proposed methods converge 
globally and characterize their convergence rates.  
Let us state the assumptions. 
\begin{Assumption} \label{ass:K}
The delayed iteration indices $\tau_i^k$ [cf.~\eqref{eq:delay}] satisfy $0 \leq k - \tau_i^k \leq K$ for all $i,k$ and for some $K \geq 0$.
\end{Assumption} 
\begin{Assumption} \label{ass:F} The function
$F(\prm)$ is $\mu$-strongly convex and $L$-smooth with $L \geq \mu > 0$.
\end{Assumption} 	
\begin{Assumption} \label{ass:fi}
For each $i$, the Hessian of the function $f_i(\prm)$ is $L_{H,i}$-Lipschitz continuous. 
\end{Assumption} 
We define $L_H \eqdef \sum_{i=1}^m L_{H,i}$ as the Lipschitz constant for the 
Hessian of the sum function $F(\prm)$. 
The first assumption above can be satisfied when Line~\ref{ciag_n:sel} of
Algorithm~\ref{alg:aciag} is implemented with either 
cyclic function selection, \ie $i_k = (k~{\rm mod}~m)+1$, or implemented with a random shuffling step
at the beginning of every epoch \cite{gurbuzbalaban2015random}. 
The second and the last assumptions 
are standard and they 
can be satisfied by a number of 
functions relevant to machine learning applications, e.g., the logistic loss function. 

Let us first present the convergence result for the {\sf CIAG} method --- 
taking the parameterization $\gamma = c / (\mu+L)$ for some $0 < c \leq 2$. 
We have:
\begin{Theorem} \label{thm:scvx}
Let Assumptions~\ref{ass:K}, \ref{ass:F} and \ref{ass:fi} hold. Consider the 
{\sf CIAG} method with its optimality gap sequence defined as
$\dst{k} \eqdef \| \prm^k - \prm^\star \|^2$ where $\prm^\star$ 
is the optimal solution to \eqref{eq:opt}. 
If the step size parameter, $c$, satisfies:
\beq \label{eq:step}
\begin{split} c < \min\Big\{ & 2, ~\frac{1}{K} \sqrt{ \frac{\mu L (\mu+L) }{2 L_H (L^2 (\dst{1})^{\frac{1}{2}} +4 L_H^2 (\dst{1})^{\frac{3}{2}} )} }, \Big( \frac{1}{K^4} \frac{ \mu L (\mu+L)^4 }{ 2 L_H^2 (L^4 \dst{1} + 16 L_H^4 (\dst{1})^3 )  } \Big)^{1/5}  \Big\} \eqs,
\end{split}
\eeq 
then there exists $\delta \in [1 - 2 \gamma \frac{ \mu L }{ L + \mu },1)$ such that the sequence $\{ \dst{k} \}_{k \geq 1}$ 
converges linearly,
\beq \label{eq:linear1}
\dst{k} \leq \delta^{\lceil (k-1)/(2K+1) \rceil} \dst{1} ,~\forall~ k \geq 1 \eqs.
\eeq
Moreover, there exists an upper bound sequence $\{ \bar{V}^{(k)} \}_{k \geq 1}$, which satisfies 
$\bar{V}^{(k)} \geq \dst{k}$ for all $k \geq 1$ and $\bar{V}^{(1)} = \dst{1}$ such that 
\beq \label{eq:rate}
\begin{split}
& \lim_{k \rightarrow \infty} \frac{ \bar{V}^{(k+1)} }{ \bar{V}^{(k)} } \leq 1 - 2 \gamma \frac{ \mu L }{ L + \mu }  = 1 - {\cal O}\Big( \frac{c}{\kappa(F)} \Big)\eqs.
\end{split}
\eeq
\end{Theorem}

Next, for the {\sf A-CIAG} method, we adopt the parameters:
\beq \label{eq:param}
\alpha = \frac{1 - \sqrt{\mu \gamma}}{1 + \sqrt{\mu \gamma}},~~\gamma = \frac{c}{L} \eqs,
\eeq
such that the extarpolation factor and step size are controlled by some $ 0  < c \leq 1/2$. 
Our result is summarized as follows: 
\begin{Theorem} \label{thm:main}
Let Assumptions~\ref{ass:K}, \ref{ass:F} and \ref{ass:fi} hold. Consider the {\sf A-CIAG} method 
with its optimality gap sequence defined as $\gap{k} \eqdef F( \prm^k ) - F(\prm^\star)$. 
If the step size parameter, $c$, satisfies:
\beq \label{eq:stepsize}
c < \min \left\{ \bar{c}_1, \bar{c}_2, \bar{c}_3, 1/2 \right\} \eqs,
\eeq 
where
\begin{align}
& \bar{c}_1 \eqdef \left( \frac{\sqrt{\mu}}{ \sqrt{18} K^2 L_H} \frac{L^2}{ \frac{20L^2}{\mu} (2 \gap{1})^{\frac{1}{2}} +  \big( \frac{40L_H}{\mu} \big)^2 (2 \gap{1})^{\frac{3}{2}} }  \right)^{\frac{1}{2}} \label{eq:c_aciag} ,\\
& \bar{c}_2 \eqdef \left( \frac{2\mu}{ 81 K^4 L_H^2} 
\frac{L^4}{ \big(\frac{20L^2}{\mu}\big)^2 (2 \bgap{1} ) +  \big( \frac{40L_H}{\mu} \big)^4 (2 \gap{1})^3 } \right)^{\frac{1}{4}},~ \bar{c}_3 \eqdef 
\frac{L}{\sqrt{324}  K^2 L_H \frac{(\gap{1})^{\frac{1}{2}}}{\sqrt{\mu}}
+ 1296 K^4 L_H^2 \frac{\gap{1}}{\mu^2} + \mu} , \notag
\end{align}
then the optimality gap $\gap{k}$ satisfies
\beq \label{eq:global}
\gap{k} \leq {\delta^{\lceil \frac{k}{2K+1} \rceil }} (2 \gap{1}),~\forall~k \geq 1 \eqs,
\eeq
for some $\delta < 1$; moreover, there exists an upper bound sequence
$\{ \bgap{k} \}_{k \geq 1}$ 
such that $\bgap{k} \geq \gap{k}$ for all $k \geq 1$ and $\bgap{1} = \gap{1}$, such that
\beq \label{eq:linear}
\lim_{k \rightarrow \infty} \frac{ \bgap{k+1} }{ \bgap{k} } = 1 - \sqrt{\mu \gamma} = 1 - \sqrt{ \frac{c}{ \kappa(F) } } \eqs.
\eeq
\end{Theorem}
The above theorems reveal that there are two phases of convergence for the {\sf CIAG} and
{\sf A-CIAG} methods -- one that converges at a \emph{slow} linear rate 
[cf.~\eqref{eq:linear1} and \eqref{eq:global}], and the asymptotic phase
where the algorithms converge linearly at a \emph{fast} rate 
comparable to the {\sf FG} and {\sf AFG} methods [cf.~\eqref{eq:rate} and \eqref{eq:linear}]. 
Such behavior is similar to  
the linear and superlinear convergence phases in the
Newton's method \cite{bertsekas1999nonlinear}, and 
they can be anticipated as the {\sf CIAG} and {\sf A-CIAG}
methods make use of the second order information. 

Due to the strong convexity of $F(\prm)$ [cf.~Assumption~\ref{ass:F}], the 
potential function used in the two theorems above are comparable, since 
\beq
\frac{\mu}{2} \dst{k} \leq \gap{k} \leq \frac{L}{2} \dst{k} \Longrightarrow \gap{k} = \Theta( \dst{k} ) \eqs.
\eeq
Both theorems require the step size parameters be chosen 
according to the initial conditions [cf.~\eqref{eq:step} and \eqref{eq:stepsize}]. 
The allowable range of step size $\gamma$ is inversely proportional to
$L_H$ and $\dst{1}$.
The former term $L_H$ measures the `quadratic-ness' of the component 
functions such that $L_H = 0$ if $f_i ( \prm )$ are quadratic.
The latter term $\dst{1}$ is the initial optimality gap for the algorithms, which is originated from the use of curvature information.
The favorable cases are when $L_H \approx 0$ or 
$\dst{1} \approx 0$  
such that we can take $c \approx 1$.
The latter implies
the following asymptotic convergence rates:
\beq \label{eq:asym_rate}
\begin{split}
& {\sf CIAG}: \lim_{k \rightarrow \infty} \frac{ \bar{V}^{(k+1)} }{ \bar{V}^{(k)} } = 1 - {\cal O} \Big( \frac{1}{\kappa(F)} \Big),\\
& \text{\sf A-CIAG}: \lim_{k \rightarrow \infty} \frac{ \bgap{k+1} }{ \bgap{k} } 
= 1 - \sqrt{\frac{1/2}{\kappa(F)}} \eqs.
\end{split}
\eeq
They coincide with the convergence rates achieved by the 
{\sf FG} and {\sf AFG} methods, respectively.
Since the per-iteration complexity of {\sf CIAG}, 
{\sf A-CIAG} are ${\cal O}(d^2)$,
while the complexity of {\sf FG}, {\sf AFG} are ${\cal O}(m d)$, 
the advantage of the proposed
methods is significant when $m \gg d$. 
It is also interesting to compare the convergence criterion for
the {\sf CIAG} and {\sf A-CIAG} methods. Most noticeably, the {\sf A-CIAG} criterion \eqref{eq:c_aciag} differs from
the {\sf CIAG} criterion \eqref{eq:step} for the region
specified by $\tilde{c}_3$ as the {\sf A-CIAG} requires the step size 
parameter be chosen as $c = {\cal O}(1/K^4)$ for large $K$. 
Since $K = \Theta(m)$, the initial optimality gap 
has to be as small as $\gap{1} = {\cal O}(1/m^4)$ to attain the 
fast rate as \eqref{eq:asym_rate}.
In comparison, the {\sf CIAG} method only requires $\dst{1} = {\cal O}(1/m^2)$
to attain \eqref{eq:asym_rate} due to the milder requirements in \eqref{eq:step}.  

When $\dst{1} \not\approx 0$ and $L_H \not\approx 0$, we can choose a small step size parameter $c$ at the beginning of {\sf CIAG} (or {\sf A-CIAG}) method that satisfies \eqref{eq:step} (or \eqref{eq:c_aciag}). By \eqref{eq:linear1} (or \eqref{eq:global}), for some finite $T_0$, the method is guaranteed to find a solution with $\dst{T_0+1} < \dst{1}$ after $T_0$ iterations. As such, we can `re-start' the method using $\prm^{(T_0+1)}$ as the initial point with a larger step size parameter $c$. The procedure can then be repeated until $c$ is increased to $c=2$ (or $c=\frac{1}{2}$ for {\sf A-CIAG}). In this way, the ideal asymptotic convergence rates in \eqref{eq:asym_rate} are achieved. 
Alternatively, we can employ other incremental optimization methods such as SAG \cite{schmidt2017minimizing} to `warm-start' the {\sf CIAG} or {\sf A-CIAG} methods with an initial point satisfying $\dst{1} \approx 0$. 
We remark that analyzing the exact complexity of a `re-starting' schedule is non-trivial which is beyond the scope of the current paper.

Nevertheless, from our numerical experiments in Section~\ref{sec:num}, we find that 
the practical step sizes can often be chosen more aggressively. For a wide range of problems, using a constant step size parameter $c$, the {\sf CIAG} and {\sf A-CIAG} methods converge quickly in terms of the number of iterations used and CPU time, without the aforementioned re-starting or warm-starting techniques.

\subsection{Proofs of Theorem~\ref{thm:scvx} and \ref{thm:main}}
The analysis for {\sf CIAG} and {\sf A-CIAG} methods consists 
of three steps.\vspace{.1cm}
\begin{itemize}
\item First, we carry out a perturbation analysis on the {\sf FG} and 
{\sf AFG} methods. Effectively, we view the {\sf CIAG} 
(\resp {\sf A-CIAG}) method 
as a perturbed {\sf FG} (\resp
{\sf AFG}) method with inexact gradients, where the errors are defined as:
\beq  
\begin{split}
& \textstyle {\bm e}_{\sf CIAG}^k \eqdef {\bm g}_{\sf CIAG}^k - \sum_{i=1}^m \grd f_i ( \prm^k ),~ \textstyle {\bm e}_{\sf ACIAG}^k \eqdef {\bm g}_{\sf ACIAG}^k - \sum_{i=1}^m \grd f_i ( \eprm^k ) \eqs.
\end{split}
\eeq
\item Second, we analyze an upper bound on $\| {\bm e}_{\sf CIAG}^k \|$
or $\| {\bm e}_{\sf ACIAG}^k\|$ 
in terms of 
\beq
\dst{k} \eqdef \| \prm^k - \prm^\star \|^2,~~\gap{k} \eqdef F(\prm^k) - F(\prm^\star) \eqs,
\eeq
where we shall use $\dst{k}$, $\gap{k}$ as the potential functions for 
{\sf CIAG} method and {\sf A-CIAG} method, respectively.
Here, we exploited Assumption~\ref{ass:fi} on the Lipschitz continuity of 
Hessian for the component functions.\vspace{.1cm}  
\item Third, using the bounds derived in the previous steps we study a nonlinear inequality system to derive 
the convergence criterion of $\dst{k}$ and $\gap{k}$.
The nonlinear inequality exhibits the desired $R$-linear convergence  
when $k \rightarrow \infty$.\vspace{.1cm} 
\end{itemize}
Overall, the key to our proof is to analyze the dynamics of optimality gap $\dst{k}$ 
(or $\gap{k}$) 
as a nonlinear system, whose convergence can be guaranteed by 
an appropriately chosen step size $\gamma$ and the asymptotic convergence rate
will only depend on the linear terms in the optimality gaps. 

In the following, we assume that the {\sf CIAG}, {\sf ACIAG}
methods are both initialized such that \eqref{eq:bkhk} holds. This can be done 
`on-the-fly' with a self-initialization step in \eqref{eq:ciag_n_inc_init} 
and the analysis below will hold for all $k \geq K$, \ie after a complete
pass through the dataset.\vspace{.2cm}


\textbf{Step 1}.
The first step is to analyze the {\sf CIAG} (\resp {\sf A-CIAG}) method 
as a perturbed version of the {\sf FG} (\resp {\sf AFG}) method, which employs 
${\bm g}_{\sf CIAG}^k$ (\resp ${\bm g}_{\sf ACIAG}^k$)
as the gradient surrogate. We have:
\begin{Prop} \label{prop:pb_ciag} Consider the {\sf CIAG} method. Under Assumption~\ref{ass:F},
if $\gamma \leq 2/(\mu + L)$, we have that for all $k \geq 1$,
\beq \label{eq:step1_ciag}
\begin{split}
\dst{k+1} & \leq
\Big( 1 - 2 \gamma \frac{\mu L }{\mu + L} \Big) \dst{k} + \gamma^2 \| {\bm e}_{\sf CIAG}^k \|^2 + 2 \gamma \sqrt{\dst{k}} \| {\bm e}_{\sf CIAG}^k \| \eqs.
\end{split}
\eeq
\end{Prop}
The proof largely follows from \cite[Section 3.3]{gurbuzbalaban2017convergence}
and is omitted.
\begin{Prop} \label{prop:pb_aciag} Consider the {\sf A-CIAG} method. Under Assumption~\ref{ass:F},
if $\gamma \leq 1/(2L)$, we have that for all $k \geq 1$,
\beq \label{eq:step1}
\begin{split}
\gap{k+1} & \leq 2 (1 - \sqrt{\mu \gamma})^k \gap{1} + \sum_{\ell=1}^k (1 - \sqrt{\mu \gamma})^{k-\ell} \Big( \sqrt{2 \gamma \gap{\ell} } \| {\bm e}_{\sf ACIAG}^\ell \|  + \sqrt{\frac{9\gamma}{\mu}} \| {\bm e}_{\sf ACIAG}^\ell \|^2
- \frac{\mu}{4} \frac{1-\mu \gamma}{\sqrt{\mu \gamma}} \| \prm^\ell - \eprm^\ell \|^2 \Big) .
\end{split}
\eeq
\end{Prop}
From the propositions, we observe that when the gradient errors vanish 
with ${\bm e}^\ell = {\bm 0}$, 
setting $\gamma = 2/(\mu+L)$ (\resp $\gamma = 1/(2L)$) gives 
$\dst{k+1} = {\cal O}( (1-1/\kappa(F))^k )$ 
(\resp $\gap{k+1} = {\cal O}( (1-\sqrt{1/2\kappa(F)})^k )$).
In other words, the linear convergence rates for {\sf FG} and {\sf AFG} methods
can be recovered.

Comparing the two propositions reveals the differences in error structure
between the {\sf CIAG} and {\sf A-CIAG} methods. First,   in \eqref{eq:step1} 
 the {\sf A-CIAG}'s error is 
convolved with the linearly converging sequence
$(1-\sqrt{\mu \gamma})^\ell$, while in \eqref{eq:step1_ciag} 
the {\sf CIAG}'s error is simply additive;
second,  \eqref{eq:step1} consists of a negative term depending on 
$\| \prm^\ell - \eprm^\ell \|^2$. This negative term revealed through our refined analysis of the error dynamics of {\sf A-CIAG} [see \eqref{eq:lowerbdphi}] and the resultant proposition is an improvement over \cite{schmidt2011convergence}.
Significantly, this negative term is an important key for establishing the fast convergence of {\sf A-CIAG}.\vspace{.1cm}

\textbf{Step 2}.
Our next step relates the gradient errors
${\bm e}_{\sf CIAG}^\ell$, 
${\bm e}_{\sf ACIAG}^\ell$
to the optimality gaps $\dst{\ell}$, $\gap{\ell}$ for the {\sf CIAG},
{\sf A-CIAG} methods, respectively. 
We obtain bounds for the errors as follows:
\begin{Prop} \label{prop:err_ciag}
Consider the {\sf CIAG} method. Under Assumptions~\ref{ass:K}, \ref{ass:F} and \ref{ass:fi}, we have that for all $\ell \geq 1$,
\beq \label{eq:bdek_ciag} \begin{split}
\| {\bm e}_{\sf CIAG}^\ell \| & \leq \gamma^2 K^2 L_H \Big( L^2 \hspace{-.3cm} \max_{ (\ell-K)_{++} \leq q \leq \ell-1 } \hspace{-.2cm} \dst{q} + 4 L_H^2 \hspace{-.1cm}
\max_{ (\ell-2K)_{++} \leq q \leq \ell-1 } \hspace{-.1cm} (\dst{q})^2 \Big) .
\end{split}
\eeq
\end{Prop}
\begin{Prop} \label{prop:err}
Consider the {\sf A-CIAG} method. Under Assumptions~\ref{ass:K}, \ref{ass:F} and \ref{ass:fi}, we have that for all $\ell \geq 1$,\vspace{-.2cm}
\beq \label{eq:bdek} \begin{split}
\| {\bm e}_{\sf ACIAG}^\ell \| & \leq \frac{3 K L_H}{2} \hspace{-.3cm} \sum_{j = (\ell-K)_{++}}^{\ell-1}  \hspace{-.3cm} \|\prm^{j+1} - \eprm^{j+1} \|^2 \\
& \hspace{-1cm} + \gamma^2 \frac{3 K^2 L_H}{2} \frac{20L^2}{\mu} \hspace{-.3cm}
\max_{ (\ell-K-1)_{++} \leq q \leq \ell-1} \hspace{-.3cm} \gap{q} + \gamma^2 \frac{3 K^2 L_H}{2} \Big( \frac{40L_H}{\mu} \Big)^2 \max_{ (\ell-2K-1)_{++} \leq q \leq \ell-1} (\gap{q})^2 \eqs.
\end{split}\vspace{-.2cm}
\eeq
\end{Prop}
We observe that the upper bounds for 
$\|{\bm e}_{\sf CIAG}^k\|$ and $\| {\bm e}_{\sf ACIAG}^k\|$ obey similar structure
in terms of
their dependences on $\dst{q}$ and $\gap{q}$. 
The upper bound for
$\| {\bm e}_{\sf ACIAG}^\ell \|$ depends on $\{ \| \prm^{j+1} - \eprm^{j+1} \|^2 \}_{j=(\ell-K)_{++}}^{\ell-1}$ which is 
the difference between the extrapolated variables and the main variables.\vspace{.1cm}

\textbf{Step 3}. We  omit the constants 
which are obtained from the propositions given in the previous 
steps. Instead, we focus on the main steps in the analysis and relegate
the exact analysis to Appendices~\ref{app:step3_ciag} and \ref{app:step3_aciag}.  

For the {\sf CIAG} method, simply substituting \eqref{eq:bdek_ciag} into \eqref{eq:step1_ciag} 
yields:\vspace{-.2cm}
\beq \label{eq:preneg_ciag}
\begin{split}
& \dst{k+1} \leq \Big( 1 - 2 \gamma \frac{\mu L }{\mu + L} \Big) \dst{k} + {\cal O}( \gamma^3 )  \max_{ (k-2K)_{++} \leq q \leq k } 
\big( (\dst{q})^{\frac{3}{2}} + (\dst{q})^2 + (\dst{q})^{\frac{5}{2}} + (\dst{q})^4 \big) \eqs,
\end{split}\vspace{-.2cm}
\eeq
where the exact form of the system will be shown in \eqref{eq:exact_ciag}.
Define  \vspace{-.1cm}
\beq \textstyle
\dst{k}_{\rm max} \eqdef \max_{ (k-2K)_{++} \leq q \leq k } \dst{q} \vspace{-.1cm}
\eeq
and using the fact that, when $\dst{k}$ is small, the second part in the 
right hand side of \eqref{eq:preneg_ciag} can be bounded 
by its lowest order term ${\cal O}(\gamma^3) (\dst{k}_{\rm max})^{\frac{3}{2}}$,
we have:
\beq \label{eq:ciag_fin}
\dst{k+1} \leq \Big( 1 - 2 \gamma \frac{\mu L }{\mu + L} \Big) \dst{k}
+ {\cal O}(\gamma^3) (\dst{k}_{\rm max})^{\frac{3}{2}} \eqs.
\eeq 

On the other hand, for the {\sf A-CIAG} method,
substituting \eqref{eq:bdek} into \eqref{eq:step1} 
and rearranging terms show 
that the $(k+1)$th optimality gap is bounded as:\vspace{-.2cm}
\beq \label{eq:preneg}
\begin{split}
& \gap{k+1} \leq 2 (1 - \sqrt{\mu \gamma})^k \gap{1} + \sum_{\ell=1}^k (1 - \sqrt{\mu \gamma})^{k-\ell} \times \\
&  \Bigg\{ {\cal O}( \gamma^{\frac{5}{2}} ) \max_{ (\ell-2K-1)_{++} \leq q \leq \ell } 
\big( (\gap{q})^{\frac{3}{2}} + (\gap{q})^2 + (\gap{q})^{\frac{5}{2}} + (\gap{q})^4 \big) 
\\
& 
+ \big( \max_{ \ell \leq q \leq \min\{ \ell+K, k\} } {\cal O}( \sqrt{\gamma \gap{q} } ) - \frac{\mu}{4} \frac{1-\mu \gamma}{\sqrt{\mu \gamma}} \big) \| \prm^\ell - \eprm^\ell \|^2 \Bigg\} \eqs,
\end{split}
\eeq
whose exact form can be found in \eqref{eq:exact_aciag}.
We emphasize that the last term on the right hand side depends on 
the difference between $\prm^\ell$ and $\eprm^\ell$, which is unique for 
{\sf A-CIAG} method due to the use of extrapolated iterates. 

To finish the proof, we identify that in both \eqref{eq:ciag_fin} and 
\eqref{eq:preneg},
the
potential functions for {\sf CIAG} and {\sf A-CIAG} methods, \ie $\dst{k+1}$ and $\gap{k+1}$, 
are upper bounded by --- a constant factor ($<1$) 
multiplied by the previous potential function's value; and a \emph{high-order}
term that depends on the delayed version of the potential function's value.
With a sufficiently small step size, the effects from the later term 
vanishes as $k \rightarrow \infty$ 
and the proposed methods converge linearly at the desired rates. 

Particularly, in the case of {\sf CIAG}, we consider
the \emph{non-linear} inequality:
\beq \label{eq:ineq2} 
R^{(k+1)} \leq p R^{(k)} + \sum_{j=1}^J q_j \max_{ k' \in [(k-M+1)_{++},k] } (R^{(k')})^{\eta_j},~\forall~k \geq 1 \eqs,
\eeq
where $0\leq p<1$, $q_j \geq 0$, $\eta_j > 1$ for all $j$
with some $J, M < \infty$. 
We have:
\begin{Prop} \label{prop:gen_ciag}
Consider \eqref{eq:ineq2}. For some $p \leq \delta < 1$, if 
\beq \label{eq:cond2} \textstyle
p + \sum_{j=1}^J q_j (R^{(1)})^{\eta_j - 1} \leq \delta < 1 \eqs,
\eeq 
then (a) $\{R^{(k)} \}_{k \geq 1}$ converges linearly as
$R^{(k)} \leq \delta^{ \lceil k / M \rceil } R^{(1)}$ for all $k \geq 1$; 
(b) there exists an upper bound sequence $\{ \bar{R}^{(k)} \}_{k \geq 1}$ with 
$\bar{R}^{(k)} \geq R^{(k)}$ for all $k \geq 1$ and $\bar{R}^{(1)} = R^{(1)}$, that converges
linearly at rate $p$ asymptotically, 
\beq \label{eq:lem2}
\begin{split}
& \lim_{k \rightarrow \infty} \bar{R}^{(k+1)} / \bar{R}^{(k)} = p \eqs.
\end{split}
\eeq
\end{Prop}
Consequently, Theorem~\ref{thm:scvx} can be proven by identifying
that $R^{(k)} = \dst{k}$ and substituting
the appropriate constants in Proposition~\ref{prop:gen_ciag}, see Appendix~\ref{app:step3_ciag}.

In the case of {\sf A-CIAG}, 
it can be verified that under our choice of step size, 
the coefficient 
in front of $\| \prm^\ell - \eprm^\ell \|^2$ is always negative for all $\ell \geq 1$.
Define
\beq
\gap{\ell}_{\rm max} \eqdef \max_{ (\ell-2K-1)_{++} \leq q \leq \ell } \gap{q} \eqs.
\eeq
When $\gap{\ell}_{\rm max}$ is small, the terms inside the last bracket of
the summation in \eqref{eq:preneg} can be bounded by its lowest order term as 
${\cal O}( \gamma^{\frac{5}{2}}  ) (\gap{\ell}_{\rm max})^{\frac{3}{2}}$. 
Therefore, we consider an upper bound sequence 
$\{ \bgap{k} \}_{k \geq 1}$ 
with
$\bgap{k} \geq \gap{k}$ for all $k$ and:
\beq \begin{split}
\bgap{k+1} & = 2( 1 - \sqrt{\mu \gamma} )^k \bgap{1}   + \sum_{\ell=1}^k (1 - \sqrt{\mu \gamma})^{k-\ell} {\cal O}( \gamma^{\frac{5}{2}} ) (\bgap{\ell}_{\rm max})^{\frac{3}{2}} \eqs.
\end{split}
\eeq 
Subtracting $(1-\sqrt{\mu \gamma}) \bar{h}^{(k)}$ from both sides gives:
\beq \label{eq:aciag_fin}
\bar{h}^{(k+1)}  = (1-\sqrt{\mu \gamma}) \bar{h}^{(k)} + {\cal O}( \gamma^{\frac{5}{2}} ) (\bgap{k}_{\rm max})^{\frac{3}{2}} \eqs,
\eeq
which resembles \eqref{eq:ciag_fin} in the case of {\sf CIAG}. 
Similar to the previous developments, we expect the system to converge
linearly at the rate $1 - \sqrt{\mu \gamma}$. 

Formally, 
consider an abstracted form of \eqref{eq:preneg}
with 
the non-negative sequence $\{ \fun{k} \}_{k \geq 1}$ that satisfies:
\beq \label{eq:abstract_sys}
\begin{split}
\fun{k+1} \leq & ~p^k b \fun{1} +  \sum_{\ell=1}^k p^{k-\ell} \Bigg\{ \sum_{j=1}^J s_j \hspace{-.1cm} \max_{ (\ell- M)_{++} \leq q \leq \ell } (\fun{q})^{\eta_j} \\
& \hspace{3.5cm} + \big( \max_{\ell \leq q \leq k} f(\fun{q}) - \bar{f}   \big) D^{(\ell)} 
\Bigg\},
\end{split}
\eeq
for all $k \geq 1$, 
where $f( \fun{q} )$ is a non-decreasing function of $\fun{q}$
and $\eta_j > 1$ for all $j$. 
The parameters 
$p, s_j, \bar{f}, b$ are all non-negative, we have $b \geq 1$ and 
$M  < \infty$, and $\{ D^{(\ell)} \}_{\ell \geq 1}$ is an arbitrary 
non-negative sequence.
The above system converges 
linearly at a rate given by the constant factor $p < 1$:
\begin{Prop} \label{prop:gen_aciag}
Consider \eqref{eq:abstract_sys}. Suppose that 
\beq \label{eq:condition} \textstyle
\bar{f} \geq f( b \fun{1} )~~~\text{and}~~~\delta \eqdef p + \sum_{j=1}^J s_j (b \fun{1})^{\eta_j - 1} < 1 \eqs.
\eeq	
Then, there exists an upper bound
sequence $\{ \bfun{k} \}_{k \geq 1}$ satisfying
\beq \label{eq:statement_aciag}
\begin{split}
& \text{(i)}~\bfun{k} \geq \fun{k},~\forall~k \geq 1,~~
\text{(ii)}~\bfun{k+1} \leq \delta^{\lceil k / M\rceil} (b \bfun{1}),~\forall~k \geq 1 \\
& \text{and}~~\text{(iii)}~
\lim_{ k \rightarrow \infty }  \bfun{k+1} / \bfun{k} = p \eqs.
\end{split}
\eeq
\end{Prop}
Finally, Theorem~\ref{thm:main} can be proven by identifying $\fun{k} = \gap{k}$ 
and substituting the appropriate constants, see Appendix~\ref{app:step3_aciag}.

\section{Numerical Experiments} \label{sec:num}
This section covers the performance of {\sf CIAG} and {\sf A-CIAG} 
using  numerical experiments.
We focus on the logistic regression problem for training  
linear classifiers. We are given $m$ data tuples $\{ ({\bm x}_i, y_i) \}_{i=1}^m$, 
where ${\bm x}_i \in \RR^d$ is the feature vector and $y_i \in \{ \pm 1 \}$ is the label.
The $i$th component function is:
\beq 
f_i (\prm) = \frac{1}{2m} \| \prm \|^2 + \log ( 1 + {\rm exp} ( -y_i \langle \prm, {\bm x}_i \rangle ) ) \eqs.
\eeq
This function has the form of a linear model in \eqref{eq:linearmodel} and it  
satisfies Assumptions~\ref{ass:F} and \ref{ass:fi}. 
For instance, an upper bound to the gradient and Hessian 
smoothness 
of $F(\prm)$ and $f_i(\prm)$, respectively, 
can be evaluated as:
\beq 
L = 1 + \frac{1}{4} \sum_{i=1}^m \| {\bm x}_i \|_2^2,~~L_{H,i} =\frac{1}{4} \| {\bm x}_i \|_2^2 \eqs.
\eeq 


\begin{figure}[t]
\centering 
\includegraphics[width=.425\linewidth]{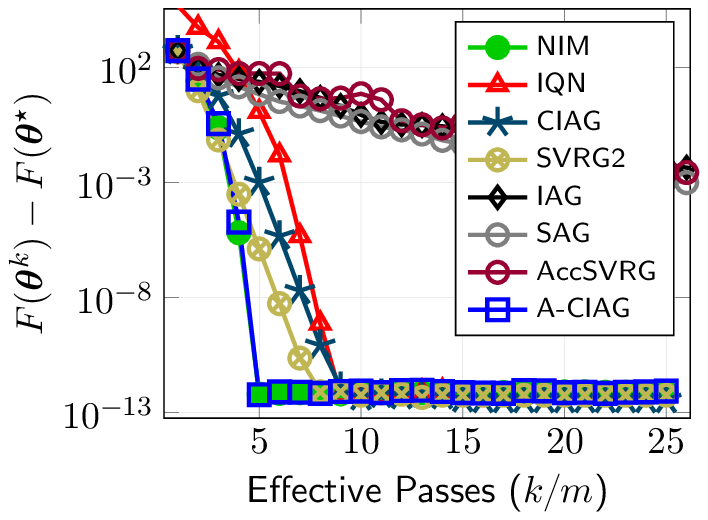}~~
\includegraphics[width=.425\linewidth]{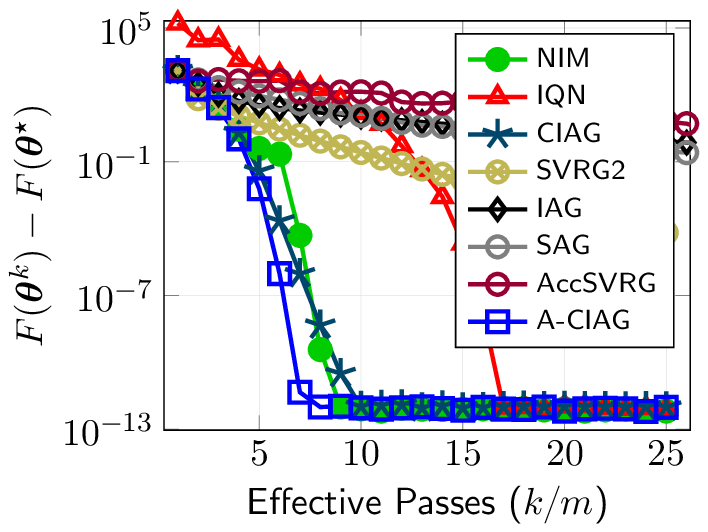}
\caption{The $y$-axis denotes the optimality
gap plotted in log-scale and the $x$-axis shows the number of effective passes
(defined as $k/m$). 
(Left) for the case with $m=1000$, $d=51$. 
(Right) for the case with $m=2000$, $d=501$.
Here, for {\sf NIM} and {\sf SVRG2} we use the step sizes $\gamma = 0.01$ and $\gamma = 0.1/L$, respectively.} \label{fig:small}\vspace{-.2cm}
\end{figure}

\subsection{Synthetic Data} 
We adopt a simple random data model. 
First, we generate $\prm_{\sf true} \sim {\cal U}[-1,1]^d$  
and the feature vector 
as ${\bm x}_i = [ \tilde{\bm x}_i;~1 ]$ where 
$\tilde{\bm x}_i \sim {\cal U}[-1,1]^{d-1}$; then, the label  is computed as
$y_i = {\rm sign}( \langle {\bm x}_i, \prm_{\sf true} \rangle )$. 

To set up the benchmark, 
the step sizes for {\sf NIM} and {\sf IQN} are  $\gamma = 1$. 
For the {\sf IAG} method, we set $\gamma = 50 / (mL)$. 
For the {\sf CIAG} and {\sf A-CIAG} methods, we set $\gamma = 1 / L$
and we set the extrapolation weight for {\sf A-CIAG} as $\alpha = 0.95$.
The above methods are implemented with
deterministic, cyclic component function selection, \ie $i_k = (k~{\rm mod}~m)+1$. 
We also compare a few stochastic methods: 
for {\sf SAG} and {\sf AccSVRG}, we set $\gamma = 50 / (mL)$; 
and the batch size
is  $B=5$ for {\sf AccSVRG} with an epoch length of $m$. 
For {\sf SVRG2}, we set an epoch length of $0.1m$. 
The step sizes for {\sf NIM} and {\sf SVRG2} will be specified later.

The evolution of the optimality gap against the number of effective passes 
through data
are shown in Figure \ref{fig:small}  for different problem sizes. 
We defined the number of effective passes as the number of \emph{iterations}
($k$)
divided by $m$. 
From Figure \ref{fig:small}, 
we observe that both {\sf NIM} and {\sf IQN} 
methods have the fastest convergence, since both methods are 
shown to converge superlinearly. 
However, we note that the curvature aided methods, 
{\sf A-CIAG}, {\sf CIAG} and {\sf SVRG2}, also demonstrate 
similar convergence speed in terms of the number of effective 
data passes used.
Especially, the speed of the proposed {\sf A-CIAG} almost matches that of the {\sf NIM}
method.


\begin{table}[h]
\begin{center}
\renewcommand\baselinestretch{1.25}\selectfont
{\footnotesize
\begin{tabular}{l l l l l}
\toprule
Dataset & {\sf A-CIAG}  & {\sf CIAG} & {\sf NIM}  & {\sf SAG}$^\ddagger$ 
\\
\midrule
\texttt{mushrooms} & \cellcolor{asublue!7}\bfseries $\bf 5.22$ pass & $43.5$ pass & $4.81$ ($4.92$) pass  & $359.9$ pass 
 \\
($m=8124, d = 112$) & \cellcolor{asublue!7}\bfseries $\bf 0.299$ sec. & $2.509$ sec. &  $1.01$ ($0.329$) sec. & $1.521$ sec. 
\\
\hline 
\texttt{a9a} & \cellcolor{asublue!7}\bfseries $\bf 3.6$ pass & $52.2$ pass & $3.0$ ($3.2$) pass & $165.8$ pass 
\\
($m=32561, d = 123$) & \cellcolor{asublue!7}\bfseries $\bf 1.067$ sec. & $15.26$ sec. & $3.38$ ($1.10$) sec. & $3.685$ sec.  
\\
\hline 
\texttt{SUSY} & \cellcolor{asublue!7}\bfseries $\bf 7.1$ pass & $7.6$ pass & $6.2$ ($6.2$) pass & $52.3$ pass 
 \\
($m=5 \times 10^6, d = 18$) & \cellcolor{asublue!7}\bfseries $\bf 24.88$ sec. & $26.81$ sec. & $35.92$ ($29.82$) sec. & $99.00$ sec. 
 \\
\hline 
\texttt{covtype} & \cellcolor{asublue!7}\bfseries $\bf 4.5$ pass & $13.5$ pass & $4.0$ ($4.5$) pass & $101.9$ pass 
\\
($m=581012, d = 54$) & \cellcolor{asublue!7}\bfseries $\bf 5.888$ sec. & $17.71$ sec. & $13.84$ ($7.06$) sec. & $32.33$ sec. 
\\
\hline
\texttt{w8a} & \cellcolor{asublue!7}\bfseries $\bf 5.5$ pass &  $7.2$ pass & $5.3$ ($5.4$) pass & $251.01$ pass 
\\
($m=49749, d = 300$) & \cellcolor{asublue!7}\bfseries $\bf 12.73$ sec. & $16.48$ sec. & $69.13$ ($13.82$) sec. & $23.20$ sec. 
\\
\hline
\texttt{mnist} & $4.3$ pass &  $143.6$ pass & \cellcolor{asublue!7}\bfseries $3.8$ ($\bf 3.8$) pass & $\geq 10^3$ pass 
\\
($m=60000, d = 784$) & $89.59$ sec. & $2801$ sec. & \cellcolor{asublue!7}\bfseries $755.2$ ($\bf 86.94$) sec. & $\geq 392$ sec.  
\\
\hline
\texttt{alpha} & \cellcolor{asublue!7}\bfseries $\bf 2.4$ pass & $7.6$ pass & $2.3$ ($2.5$) pass & $80.5$ pass 
\\[-.05cm]
($m=5 \times 10^5, d = 500$) & \cellcolor{asublue!7}\bfseries $\bf 149.5$ sec. & $475.6$ sec. &  $1111.6$ ($176.4$) sec. & 210.7 sec. 
\\
\bottomrule
\end{tabular}}\vspace{-.2cm}
\end{center}
\caption{Performance comparison. 
We show the number of effective passes (defined as $k/m$) 
and the wall clock time required to reach convergence 
with $\| \grd F(\prm^k) \| \leq 10^{-10}$. For the {\sf NIM} method,
we tested both `exact' and `inexact' settings in Hessian 
inverse. Results inside the brackets $(\cdot)$ correspond
to `inexact' setting. ($^\ddagger$Results for {\sf SAG} are averaged 
over $10$ trials.)}\label{tab:real} \vspace{-.2cm}
\end{table}

\begin{table}[h]
\begin{center}
\renewcommand\baselinestretch{1.25}\selectfont
\setlength\tabcolsep{4pt}
{\scriptsize
\begin{tabular}{l l l l l l l l}
\toprule
 & \texttt{mushrooms} & \texttt{a9a} & \texttt{SUSY} & \texttt{covtype} & \texttt{w8a} & \texttt{mnist} & \texttt{alpha} \\ 
\midrule
{\sf A-CIAG} &$\gamma = \frac{10^{-3}}{L}$&$\gamma = \frac{10^{-4}}{ L}$&$\gamma = \frac{10^{-5}}{ L}$&$\gamma = \frac{5 \cdot 10^{-6}}{L}$&$\gamma = \frac{10^{-3}}{L}$&$\gamma = \frac{10^{-4}}{L}$&$\gamma = \frac{5 \cdot 10^{-6}}{ L}$\\
& $\alpha = 0.99$& $\alpha = 0.99$ & $\alpha = 0.99$ & $\alpha =0.99$ & $\alpha = 0.99$
& $\alpha = 0.99$ & $\alpha = 0.99$ \\
\midrule
{\sf CIAG}&$\gamma = \frac{10^{-3}}{L}$&$\gamma = \frac{2 \cdot 10^{-4}}{ L}$&$\gamma = \frac{10^{-5}}{ L}$&$\gamma = \frac{5 \cdot 10^{-6}}{L}$&$\gamma = \frac{10^{-3}}{L}$&$\gamma = \frac{10^{-4}}{L}$&$\gamma = \frac{ 10^{-5}}{ L}$\\
\bottomrule
\end{tabular}\vspace{-.2cm}}
\end{center}
\caption{Parameters used for {\sf A-CIAG}, {\sf CIAG}
methods in the numerical experiments.} \label{tab:stepsize}\vspace{-.6cm}
\end{table}

\subsection{Real Data} We empirically compare the algorithms on the
datasets from \texttt{LibSVM} \cite{CC01a}. 
For this example, the algorithms are implemented in \texttt{C++} based on 
the source codes by \cite{rodomanov2016superlinearly} (available: \url{http://github.com/arodomanov/nim_icml16_code}) to demonstrate the fastest possible practical performance.
We only compare 
the proposed {\sf CIAG}, {\sf A-CIAG} to {\sf IAG}, {\sf NIM} and {\sf SAG},
where the first four algorithms employ a
deterministic, cyclic component function selection
and the last algorithm employs random component function selection. 
For the {\sf NIM} method, we tested both of its \emph{exact}
and \emph{inexact} setting, where the \emph{inexact} setting is a double-loop
method which uses 
a conjugate gradient method to tackle the Hessian inverse involved. 
We have used a mini-batch setting for all the tested methods
such that each $f_i(\prm)$ is composed of $B=5$ data tuples. 
The numerical experiments were conducted on a Laptop computer with 
an Intel Core i7 2.8Ghz quad-core processor and 16 gigabytes of memory.

An overview of the performance comparison   can be 
found in Table~\ref{tab:real}, while Table~\ref{tab:stepsize} 
shows the algorithms' parameter settings used
for different dataset. 
As seen, the {\sf A-CIAG} method outperformed the benchmarks
for many of the considered datasets in terms of the wall clock convergence time,
and the number of effective passes required is comparable to the best available
method, {\sf NIM}. 
An exception is the experiment on the \texttt{mnist} dataset. 
We suspect that this is due to the potentially poor condition number with 
the \texttt{mnist} dataset, as we observe that the unaccelerated methods
{\sf SAG}, {\sf IAG} and {\sf CIAG} exhibit significantly slower convergence
than for the other datasets. 

\begin{figure}[t]
\centering
\includegraphics[height=.11\textheight]{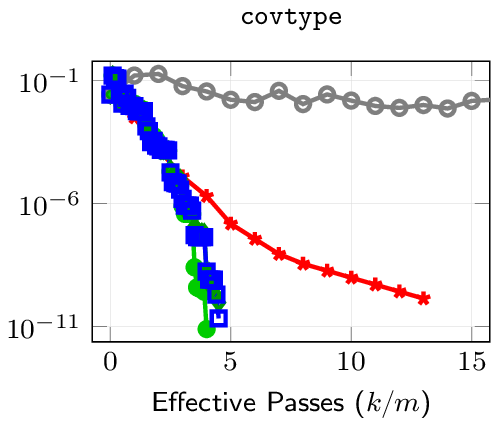}
\includegraphics[height=.11\textheight]{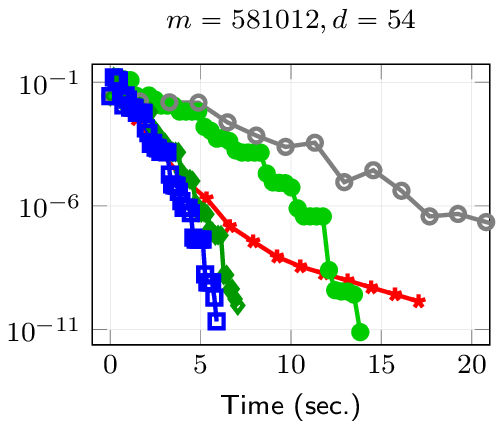}
\includegraphics[height=.11\textheight]{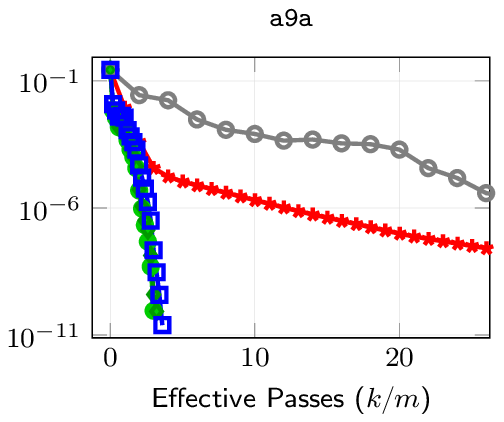}
\includegraphics[height=.11\textheight]{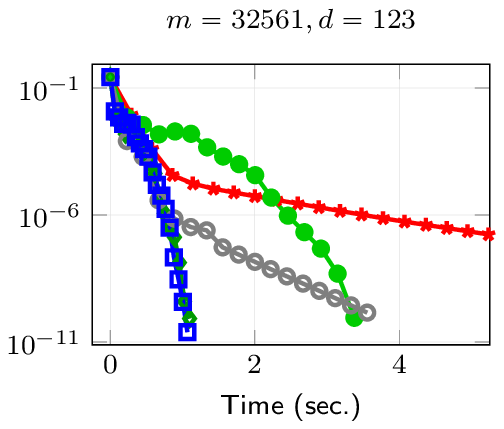}
\includegraphics[height=.11\textheight]{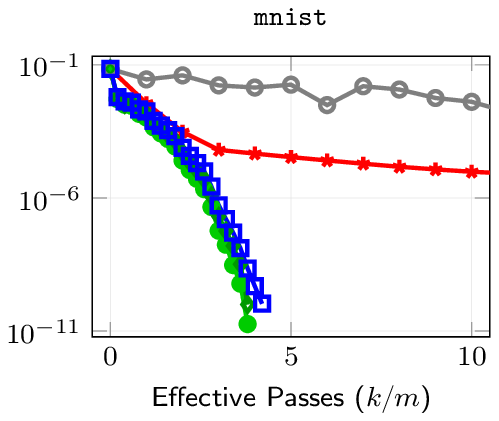}
\includegraphics[height=.11\textheight]{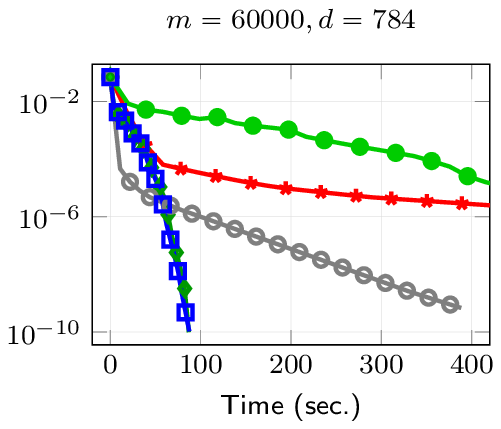}
\includegraphics[height=.11\textheight]{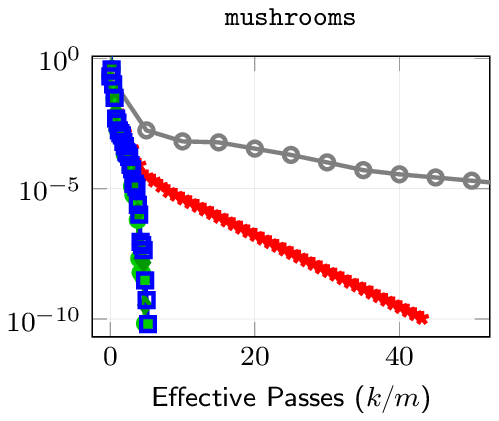}
\includegraphics[height=.11\textheight]{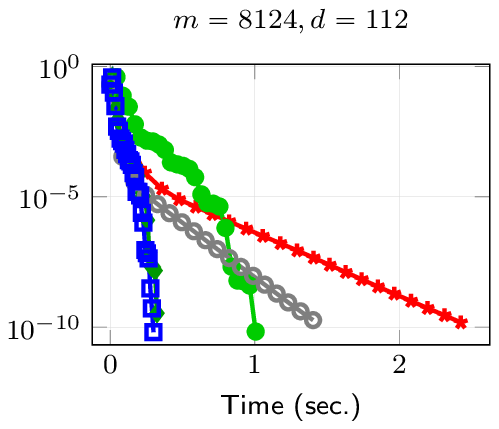}
\includegraphics[height=.11\textheight]{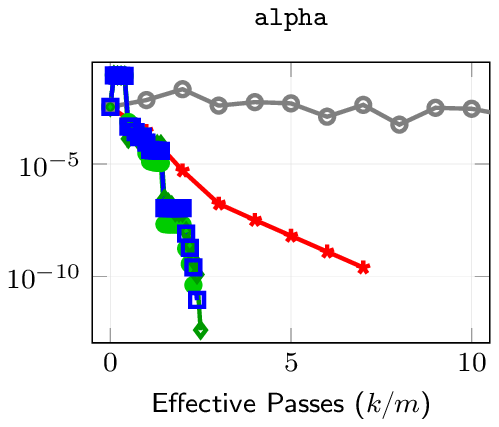}
\includegraphics[height=.11\textheight]{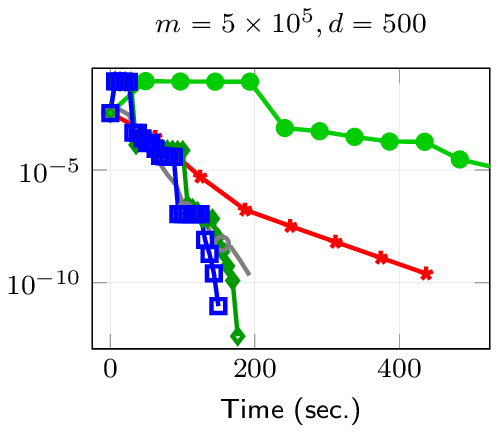}
\includegraphics[height=.11\textheight]{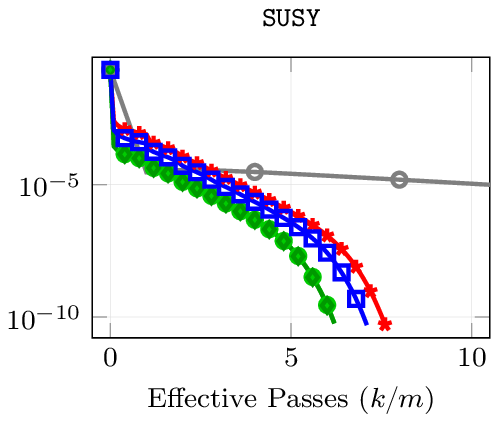}
\includegraphics[height=.11\textheight]{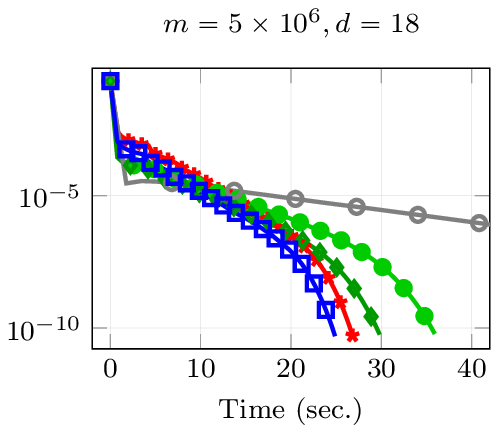}
\includegraphics[height=.11\textheight]{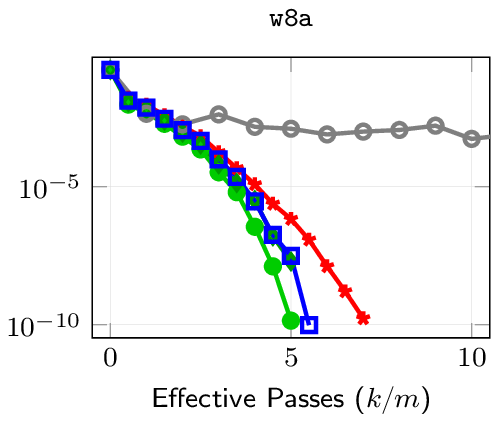}
\includegraphics[height=.11\textheight]{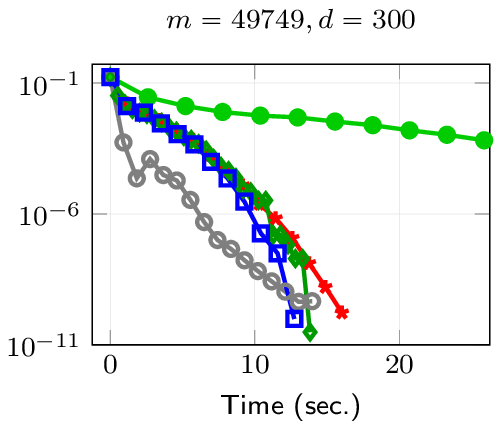}\vspace{.4cm}

\includegraphics[width=.6\linewidth]{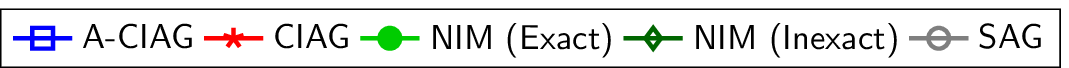}
\caption{Evolution of $\| \grd F( \prm^k ) \|$ ($y$-axis) against number of effective passes and wall clock
time on the datasets.
The experiment settings are the same as in Table~\ref{tab:real}.}\label{fig:real}\vspace{-.2cm}
\end{figure}

To investigate the behavior of the algorithms, Figure~\ref{fig:real}
shows the evolution of gradient's norm $\| \grd F( \prm^k ) \|$ 
against the number of effective passes and wall clock time
on the tested datasets.
As seen, the convergence speed of {\sf A-CIAG} matches that of the Newton-based
{\sf NIM}, yet the wall clock time required is faster
as it does not involve computing the Hessians' inverse.
It is worthwhile pointing out that except for 
the datasets \texttt{covtype} and 
\texttt{alpha}, 
the {\sf SAG} method achieves a solution accuracy of 
$\| \grd F( \prm^k ) \| \leq 10^{-6}$ in less wall clock time than
the proposed methods. 
The benefits of our proposed methods are significant when
the solution accuracy is high, as predicted by our theoretical results.\vspace{-.2cm}

\section{Conclusions}
We proposed two new optimization methods utilizing 
curvature-aided gradient tracking for large-scale optimization via incremental
data processing. 
The proposed methods, {\sf CIAG} and {\sf A-CIAG}, attain an $\epsilon$-optimal solution
with only ${\cal O}( {\kappa(F)} \log (1/\epsilon))$ and ${\cal O}( \sqrt{\kappa(F)} \log (1/\epsilon))$ \emph{iterations}, respectively,  
for a small $\epsilon$, and provided that the initial point is close to the optimal solution.
Numerical experiments on real and synthetic data demonstrate the benefit of our algorithms.
Future work includes designing and analyzing an optimal `re-starting' scheme to achieve the ideal convergence rate \eqref{eq:asym_rate} when the initial point is not close to the optimal solution, and extending the methods to tackling non-smooth optimization problems.

\bibliography{cigd_ref}
\bibliographystyle{spmpsci}

\appendix
\section{Proof of Proposition \ref{prop:err_ciag}} \label{app:errB_ciag}
Let us express the gradient error as
${\bm e}_{\sf CIAG}^k = \sum_{i=1}^m \big( \grd f_i ( \prm^{\tau_i^k} ) +
\grd^2 f_i ( \prm^{\tau_i^k} ) ( \prm^k - \prm^{\tau_i^k} ) - \grd f_i ( \prm^k ) \big)$. Applying Lemma~\ref{lem:nes}:
\beq \label{eq:tmp_eq}
\begin{split}
& \| {\bm e}_{\sf CIAG}^k \| \leq \sum_{i=1}^m \frac{L_{H,i}}{2} \| \prm^{\tau_i^k} - \prm^k \|^2 \leq \sum_{i=1}^m \frac{L_{H,i}}{2} 
\underbrace{(k - \tau_i^k)}_{\leq K} \sum_{j=\tau_i^k}^{k-1} \| \prm^{j+1} - \prm^j \|^2 \\
& \leq \frac{K L_H}{2} \sum_{j= (k-K)_{++}}^{k-1} \hspace{-.1cm}  \| \prm^{j+1} - \prm^j \|^2  \leq \frac{K L_H}{2} \gamma^2 \sum_{j=(k-K)_{++}}^{k-1} \hspace{-.1cm}  \| {\bm e}_{\sf CIAG}^j + \grd F(\prm^j) \|^2 \\
&  \leq \gamma^2 K L_H \sum_{j=(k-K)_{++}}^{k-1} \big( \| {\bm e}_{\sf CIAG}^j \|^2 + \| \grd F(\prm^j) \|^2 \big) \eqs.
\end{split}
\eeq
Furthermore, we have
\beq
\| \grd F(\prm^j) \|^2 = \| \grd F(\prm^j) - \grd F(\prm^\star) \|^2 \leq L^2 \dst{j},  
\eeq
\beq \begin{split}
\| {\bm e}_{\sf CIAG}^j \| &
 \overset{(a)}{\leq} \sum_{i=1}^m L_{H,i} \big( \dst{j} + \dst{\tau_i^j} \big)  \leq 2 L_H \max_{ \ell \in \{ \tau_i^j \}_{i=1}^m \cup \{j\} } \dst{\ell} \eqs,
\end{split}
\eeq
where (a) is due to $\| {\bm a} - {\bm b} \|^2 \leq 2 (\|{\bm a}\|^2 + \| {\bm b} \|^2)$.
Plugging these back into \eqref{eq:tmp_eq}
and using $\tau_i^{k-K} \geq k - 2K$ gives:
\beq \begin{split}
\| {\bm e}_{\sf CIAG}^k \| & \leq \gamma^2 K L_H \sum_{j=(k-K)_{++}}^{k-1} \Big( L^2 \dst{j} + 
\big( 2 L_H \max_{ \ell \in \{ \tau_i^j \}_{i=1}^m \cup \{j\} } \dst{\ell} \big)^2 \Big) \\
&  \leq \gamma^2 K^2 L_H \Big( L^2 \max_{ (k-K)_{++} \leq \ell \leq k-1 } \dst{\ell} + 4 L_H^2 
\max_{ (k-2K)_{++} \leq \ell \leq k-1 } (\dst{\ell})^2 \Big) \eqs.
\end{split}
\eeq

\section{Step 3 in the Proof of Theorem~\ref{thm:scvx}} \label{app:step3_ciag}
Combining Proposition~\ref{prop:pb_ciag} and \ref{prop:err_ciag} yields 
\beq \label{eq:exact_ciag}
\begin{split}
\dst{k+1} 
& \leq \Big( 1 - 2\gamma  \frac{ \mu L }{\mu + L}\Big) \dst{k} \\
& \hspace{0cm} + 2 \gamma^3 K^2 L_H \Big( L^2 \max_{ (k-K)_{++} \leq \ell \leq k } (\dst{\ell})^{\frac{3}{2}}+ 4 L_H^2 
\max_{ (k-2K)_{++} \leq \ell \leq k } (\dst{\ell})^{\frac{5}{2}} \Big) \\
& \hspace{0cm} + 2 \gamma^6 K^4 L_H^2 \Big( L^4 \hspace{-.1cm} \max_{ (k-K)_{++} \hspace{-.1cm} \leq \ell \leq k-1 } (\dst{\ell})^2 + 16 L_H^4 \hspace{-.1cm}
\max_{ (k-2K)_{++} \leq \ell \leq k-1 } \hspace{-.1cm} (\dst{\ell})^4 \Big),
\end{split}
\eeq
which is the exact form for Eq.~\eqref{eq:preneg_ciag}. The right hand side 
in \eqref{eq:exact_ciag}
can be decomposed into two terms --- the first term is of 
the same order as $\dst{k}$,
and the other terms are \emph{delayed} and \emph{higher-order} terms of $\dst{\ell}$. 


Observe that \eqref{eq:exact_ciag} is a special case of \eqref{eq:ineq2}
in Proposition~\ref{prop:gen_ciag} with
$R^{(k)} = \dst{k}$, $M=2K+1$, $p=1 - 2 \gamma \mu L / (\mu + L)$ and
\beq  
\begin{split}
& q_1 = 2 \gamma^3 K^2 L^2 L_H,~\eta_1 = 3/2,~q_2 = 8 \gamma^3 K^2 L_H^3,~\eta_3 = 5/2 \eqs, \\
& q_3 = 2 \gamma^6 K^4 L_H^2 L^4,~\eta_3 = 2,~q_4 = 32 \gamma^6 K^4 L_H^6,~\eta_4 = 4 \eqs.
\end{split}
\eeq
The corresponding convergence condition in \eqref{eq:cond2} can be satisfied if
\beq 
\begin{split}
& \gamma^5 \!~ 2K^4 L_H^2 \Big( L^4 \dst{1} + 16 L_H^4 (\dst{1})^3 \Big) < \frac{ \mu L } { \mu + L } \\
& \text{and}~~\gamma^2 \!~ 2K^2 L_H \Big( L^2 (\dst{1})^{1/2} + 4 L_H^2 (\dst{1})^{3/2} \Big) < \frac{ \mu L } { \mu + L } \eqs,
\end{split}
\eeq
which can be implied by \eqref{eq:step}. 
The proof is thus concluded.

\section{Proof of Proposition~\ref{prop:gen_ciag}}
The proof of the proposition is divided into two parts. 
We first show that under \eqref{eq:cond2}, the sequence 
$\{ R^{(k)} \}_{k \geq 1}$ converges linearly as in part (a) of the proposition; then we show
that the rate of convergence is asymptotically given by $p$
as in part (b) of the proposition [cf.~\eqref{eq:lem2}]. 

The first part of the proof is achieved using induction on all $\ell \geq 1$ with: 
\beq \label{eq:ind}
R^{(k)} \leq \delta^\ell \!~ R^{(1)},~\forall~k=(\ell-1)M + 2,..., \ell M + 1\eqs.
\eeq 
The base case when $\ell=1$ can be straightforwardly established:
\beq
\begin{split}
& \textstyle R^{(2)} \leq p R^{(1)} + \sum_{j=1}^J q_j (R^{(1)})^{\eta_j} \leq \delta R^{(1)} \eqs, \\
& \vdots \\
& \textstyle R^{(M+1)} \leq p R^{(M)} + \sum_{j=1}^J q_j (R^{(0)})^{\eta_j} \leq \delta R^{(1)} \eqs. 
\end{split}
\eeq
Suppose that the statement \eqref{eq:ind} is true up to $\ell=c$, for $\ell=c+1$,
we have:
\beq \notag
\begin{split}
R^{( cM+ 2)} & \leq p R^{( cM+1 )} + \sum_{j=1}^J q_j \max_{ k' \in [ (c-1)M + 2, cM +1 ] } (R^{(k')})^{\eta_j} \\
& \leq p \big( \delta^c R^{(1)} \big) + \sum_{j=1}^J q_j  \big( \delta^c R^{(1)} \big)^{\eta_j}  \leq \delta^c \!~ \Big( pR^{(1)} + \sum_{j=1}^J q_j  (R^{(1)})^{\eta_j} \Big) \leq \delta^{c+1} R^{(1)} \eqs.
\end{split}
\eeq 
Similar statement also holds for $R^{(k)}$ with $k=cM+3,...,(c+1)M+1$. We thus conclude with:
\beq
R^{(k)} \leq \delta^{ \lceil (k-1) / M \rceil } \!~ R^{(1)},~\forall~ k \geq 1 \eqs,
\eeq
which proves the first part of the proposition. 

The second part of the proof establishes the asymptotic 
linear rate of convergence in \eqref{eq:lem2}. We
consider the upper bound sequence $\{ \bar{R}^{(k)} \}_{k \geq 1}$ such that
$\bar{R}^{(1)} = R^{(1)}$ and the inequality \eqref{eq:ineq2} is tight for 
$\{ \bar{R}^{(k)} \}_{k \geq 1}$. Obviously, it also holds that
$\bar{R}^{(k)} \leq \delta^{ \lceil (k-1) / M \rceil } \bar{R}^{(1)}$
for all $k \geq 1$. 
Now, observe that
\beq \label{eq:sec_part}
\frac{\bar{R}^{(k+1)}}{\bar{R}^{(k)}} =  p + \frac{ \sum_{j=1}^J q_j \max_{ k' \in [(k-M+1)_{++}, k] } (R^{(k')})^{\eta_j} }{ \bar{R}^{(k)} } \eqs.
\eeq
For any $k' \in [k-M+1,k]$ and any $\eta > 1$, we have:
\beq \begin{split}
& \frac{ (\bar{R}^{(k')})^{\eta} }{ \bar{R}^{(k)} } = \frac{ \bar{R}^{(k')} }{ \bar{R}^{(k)} } \!~ (\bar{R}^{(k')})^{\eta-1} \leq \frac{ \bar{R}^{(k')} }{ \bar{R}^{(k)} } (R^{(1)})^{\eta-1} \delta^{ (\lceil \frac{k'-1}{M} \rceil)(\eta-1) }\eqs.
\end{split}
\eeq
Note that as $\bar{R}^{(k+1)} / \bar{R}^{(k)} \geq p$, we have:
\beq
\frac{ (\bar{R}^{(k')})^{\eta} }{ \bar{R}^{(k)} } \leq p^{-M} (R^{(1)})^{\eta-1} \delta^{ (\lceil \frac{k'-1}{M} \rceil)(\eta-1) } \eqs.
\eeq
Taking $k \rightarrow \infty$ shows that the right hand side vanishes. As a result, 
we have 
$\lim_{k \rightarrow \infty} \bar{R}^{(k+1)} / \bar{R}^{(k)} = p$. This proves 
part (b) of the proposition.

\section{Proof of Proposition \ref{prop:pb_aciag}}
The following proof is partially inspired by \cite{schmidt2011convergence,bubeck2015convex,nitanda2014stochastic}. 
For simplicity, we drop the subscript {\sf ACIAG} in 
${\bm g}_{\sf ACIAG}^k$ and ${\bm e}_{\sf ACIAG}^k$. 
Define $\rho \eqdef 1 - \sqrt{\mu \gamma}$ and the estimation sequence as:
\beq \label{eq:phi1} \begin{split} 
\Phi_1 ( \prm) & \eqdef F ( \eprm^1 ) + \frac{ \mu}{2} \| \prm - \eprm^1 \|^2 \\
\Phi_{k+1}( \prm ) & \eqdef \rho \!~\Phi_k ( \prm )   + \sqrt{\mu \gamma} \Big( 
F( \eprm^k) + \langle {\bm g}^k, \prm - \eprm^k \rangle + \frac{\mu}{2} \| 
\prm - \eprm^k \|^2 \Big) \eqs,
\end{split}
\eeq
where ${\bm g}^k \eqdef {\bm b}^k + {\bm H}^k \eprm^k$ is the 
gradient surrogate used in \eqref{eq:ciag_n_imp}. Recall that
${\bm e}^k \eqdef {\bm g}^k - \grd F( \eprm^k )$ is the gradient error. 
The following inequality, which holds for all $\prm \in \RR^d$, 
can be immediately obtained using 
\eqref{eq:phi1} and the $\mu$-strong convexity of $F(\prm)$:
\beq \label{eq:claim0}
\begin{split}
& \Phi_{k+1} (\prm) - F(\prm) =
\rho \Phi_k ( \prm ) - F(\prm) \\
& \hspace{1.5cm} + \sqrt{\mu \gamma} \Big( 
F( \eprm^k) + \langle \grd F(\eprm^k) + {\bm e}^k, \prm - \eprm^k \rangle + \frac{\mu}{2} \| 
\prm - \eprm^k \|^2 \Big)  \\
& \leq \rho \big( \Phi_k ( \prm ) - F(\prm) \big) + \sqrt{\mu \gamma}  \langle {\bm e}^k, \prm - \eprm^s \rangle \\
& \leq \rho^k \big( \Phi_1( \prm ) - F(\prm) \big) + \sum_{\ell=1}^k \rho^{k-\ell} \sqrt{ \mu \gamma } \langle {\bm e}^\ell, \prm - \eprm^\ell \rangle  \eqs.
\end{split}
\eeq
To facilitate our development, let us denote:
\beq
\Phi_k^\star \eqdef \min_{ \prm } \Phi_k ( \prm ),~~{\bm v}^k \eqdef \arg \min_{ \prm } \Phi_k ( \prm ) \eqs.
\eeq
By setting $\prm = \prm^\star$ in \eqref{eq:claim0}, we have:
\beq \label{eq:motivate}
\begin{split}
& \Phi_{k+1}^\star - F(\prm^\star) \leq \Phi_{k+1}(\prm^\star) - F(\prm^\star) \\
& \leq \rho^k \Big( \frac{\mu}{2} \| \prm^\star - \eprm^1 \|^2 
+ F(\eprm^1) - F(\prm^\star) \Big) 
+ \sum_{\ell=1}^k \rho^{k-\ell} \sqrt{ \mu \gamma } \langle {\bm e}^\ell, \prm^\star - \eprm^\ell \rangle \\
& \leq 2 \rho^k \Big( F(\prm^1) - F(\prm^\star) \Big) 
+ \sum_{\ell=1}^k \rho^{k-\ell} \sqrt{ \mu \gamma } \langle {\bm e}^\ell, \prm^\star - \eprm^\ell \rangle \eqs.
\end{split}
\eeq
Now, if $F( \prm^{k+1} ) \leq \Phi_{k+1}^\star$, then the inequality
above shows the evolution of
the optimality gap $\gap{k}$. 
This motivates our next step, relating $F( \prm^{k+1} )$ to $\Phi_{k+1}^\star$.\vspace{.2cm}

\textbf{Lower bounding $\Phi_{k+1}^\star$ in the presence of errors}. Since $\grd^2 \Phi_k ( \prm ) = \mu {\bm I}$, the function $\Phi_k(\prm)$ is quadratic 
and we can represent $\Phi_k(\prm)$ 
alternatively as
\beq \label{eq:phi2}
\Phi_k (\prm ) = \Phi_k^\star + \frac{\mu}{2} \| \prm - {\bm v}^k \|^2 \eqs.
\eeq
By 
substituting \eqref{eq:phi2} into the definition of $\Phi_{k+1} (\prm)$ in \eqref{eq:phi1}
and evaluating the first order optimality condition of the latter, 
we have:
\beq \label{eq:fonc}
\begin{split}
& \sqrt{\mu \gamma} ( {\bm g}^k + \mu ( {\bm v}^{k+1} - \eprm^k ) ) + \rho \!~ \mu ( {\bm v}^{k+1} - {\bm v}^k ) = {\bm 0} \eqs,\\
& \Longrightarrow {\bm v}^{k+1} = \rho {\bm v}^k + \sqrt{ \mu \gamma } \eprm^k - \sqrt{\frac{\gamma}{\mu}} {\bm g}^k \eqs.
\end{split}
\eeq
By setting $\prm=\eprm^k$ in \eqref{eq:phi1} 
and using the recursive definition of $\Phi_{k+1} (\prm)$, we obtain
\beq \label{eq:compare_rhs1} \begin{split}
\Phi_{k+1} ( \eprm^k ) & 
= \rho \Phi_{k} ( \eprm^k ) + \sqrt{\mu \gamma} F( \eprm^k )  
 = \rho \big( \Phi_k^\star + \frac{\mu}{2} \| \eprm^k - {\bm v}^k \|^2 \big) + \sqrt{\mu \gamma} F( \eprm^k ) \eqs,
\end{split}
\eeq
while setting $\prm=\eprm^k$ in \eqref{eq:phi2} and using \eqref{eq:fonc} gives us:
\beq \label{eq:compare_rhs2}
\begin{split}
\Phi_{k+1} ( \eprm^k )  & 
 = \Phi_{k+1}^\star + \frac{\mu}{2} \Big( \rho^2 \|
\eprm^k - {\bm v}^k \|^2 + \frac{\gamma}{\mu} \| {\bm g}^k \|^2 + 2 \rho \sqrt{\frac{\gamma}{\mu}} \langle {\bm g}^k, \eprm^k - {\bm v}^k \rangle \Big) \eqs.
\end{split}
\eeq
Comparing the right hand side of \eqref{eq:compare_rhs1} and \eqref{eq:compare_rhs2} shows:
\beq \notag \begin{split}
\Phi_{k+1}^\star & = \rho \Big( \Phi_k^\star + \frac{\mu}{2} \| \eprm^k - {\bm v}^k \|^2 \Big) + \sqrt{\mu \gamma} F( \eprm^k ) \\
& 
~~~ - \frac{\mu}{2}\Big( \rho^2 \|
\eprm^k - {\bm v}^k \|^2 + \frac{\gamma}{\mu} \| {\bm g}^k \|^2 + 2 \rho \sqrt{\frac{\gamma}{\mu}} \langle {\bm g}^k, \eprm^k - {\bm v}^k \rangle \Big) \\
& = \rho \Phi_k^\star + \sqrt{\mu \gamma} F(\eprm^k) 
+ \frac{\mu}{2} \rho \sqrt{\mu \gamma} \| \eprm^k - {\bm v}^k \|^2  - \frac{\gamma}{2} \| {\bm g}^k \|^2 - \rho \sqrt{\mu \gamma} \langle {\bm g}^k, \eprm^k - {\bm v}^k \rangle \eqs.
\end{split}
\eeq
Using the fact   
${\bm v}^k - \eprm^k = (\sqrt{\mu \gamma})^{-1} \big( \eprm^k - \prm^k \big)$
(proven in Section~\ref{app:equality}), we have
\beq \label{eq:phistar} \begin{split}
\Phi_{k+1}^\star & = \rho \Phi_k^\star + \sqrt{\mu \gamma} F(\eprm^k) + \frac{\mu}{2} \frac{ \rho }{ \sqrt{\mu \gamma} } \| \eprm^k - \prm^k \|^2  
- \frac{ \gamma}{2} \| {\bm g}^k \|^2 - \rho \langle {\bm g}^k,  \prm^k - \eprm^k \rangle \eqs.
\end{split}
\eeq
We obtain the following chain:
\beq \label{eq:lowerbdphi}
\begin{split}
& F( \prm^{k+1} ) - \Phi_{k+1}^\star \overset{(a)}{\leq} F( \eprm^k ) - \gamma \langle \grd F( \eprm^k ), {\bm g}^k \rangle + \frac{L \gamma^2}{2} \| {\bm g}^k \|^2  - \Phi_{k+1}^\star \\
& \overset{(b)}{=} \rho \!~ \Big( F(\eprm^k )  + \langle {\bm g}^k,  \prm^k - \eprm^k \rangle - \Phi_k^\star \Big) \\
& \hspace{.5cm} -\gamma \langle \grd F( \eprm^k ), {\bm g}^k \rangle 
+ \frac{\gamma}{2} \Big( 1 + L \gamma \Big) \|{\bm g}^k \|^2 - \frac{\mu}{2} \frac{  \rho }{ \sqrt{\mu \gamma} } \| \eprm^k - \prm^k \|^2 \\
& \overset{(c)}{=} \rho \!~ \big( F(\eprm^k)  + \langle \grd F(\eprm^k),  \prm^k - \eprm^k \rangle 
- \Phi_k^\star \big)  
-\gamma \langle \grd F( \eprm^k ), {\bm g}^k \rangle \\
& \hspace{.5cm} + \rho \langle {\bm e}^k, \prm^k - \eprm^k \rangle + \frac{\gamma}{2} \Big( 1 + L \gamma \Big) \|{\bm g}^k \|^2 - \frac{\mu}{2} \frac{ \rho  }{ \sqrt{\mu \gamma} } \| \eprm^k - \prm^k \|^2 \\
& \overset{(d)}{\leq} \rho \!~ \big( F(\prm^k) - \Phi_k^\star + \langle {\bm e}^k, \prm^k - \eprm^k \rangle \big) - \frac{\mu}{2} \frac{ 1 - \mu \gamma}{ \sqrt{\mu \gamma} } \| \eprm^k - \prm^k \|^2 \\
& \hspace{.5cm} + \frac{\gamma}{2} \Big( 1 + L \gamma \Big) \|{\bm g}^k \|^2 - \gamma \langle \grd F( \eprm^k ), {\bm g}^k \rangle \\
& \overset{(e)}{\leq} \rho \!~ \big( F(\prm^k) - \Phi_k^\star + \langle {\bm e}^k, \prm^k - \eprm^k \rangle \big)  - \frac{\mu}{2} \frac{ 1 - \mu \gamma}{ \sqrt{\mu \gamma} } \| \eprm^k - \prm^k \|^2   + \gamma \| {\bm e}^k \|^2 \eqs,
\end{split}
\eeq
where (a) is due to the $L$-smoothness of $F$; (b)  
is due to \eqref{eq:phistar}; 
(c) is obtained by expanding ${\bm g}^k$
as $\grd F(\eprm^k) + {\bm e}^k$; 
(d) is obtained by adding and subtracting $(\mu/2) \| \prm^k - \eprm^k \|^2$ inside
the first bracket, applying the identity $ \rho + \rho / \sqrt{\mu \gamma} = (1 - \mu \gamma) / \sqrt{\mu \gamma}$, and using the $\mu$-strong convexity
of $F$; 
and (e) is due to the following chain of inequalities:
\beq \notag \begin{split}
& \frac{\gamma}{2} \Big( 1 + L \gamma \Big) \|{\bm g}^k \|^2 - \gamma \langle \grd F( \eprm^k ), {\bm g}^s \rangle \\
& \leq \frac{\gamma}{2} \Big( 1 + L \gamma \Big) \Big( \| {\bm e}^k \|^2 + \| \grd F( \eprm^k ) \|^2  \Big) + \frac{ L \gamma^2 }{2} \Big( \| \grd F(\eprm^k ) \|^2 + \| {\bm e}^k \|^2 \Big)  - \gamma \| \grd F( \eprm^k ) \|^2 \\
& = \Big( \frac{\gamma}{2} + L \gamma^2 \Big) \| {\bm e}^k \|^2 +
 \Big( -\frac{\gamma}{2} + L \gamma^2 \Big) \| \grd F( \eprm^k ) \|^2  
\leq \gamma \| {\bm e}^k \|^2 \eqs.
\end{split}
\eeq
As $\Phi_1( \prm^1 )  = F( \prm^1 ) = \Phi_1^\star$, applying the inequality 
\eqref{eq:lowerbdphi} recursively shows:
\beq \label{eq:finlowerbd}
\begin{split}
& F( \prm^{k+1} ) - \Phi_{k+1}^\star \leq \\
& 
\sum_{\ell=1}^k \rho^{k-\ell} 
\Big( (1-\sqrt{\mu \gamma})
\langle {\bm e}^\ell, \prm^\ell - \eprm^\ell \rangle 
+ \gamma \|{\bm e}^\ell \|^2 - \frac{\mu}{2} \frac{1-\mu \gamma}{\sqrt{\mu \gamma}} \| \eprm^\ell - \prm^\ell \|^2 \Big) \eqs.
\end{split}
\eeq
Importantly, \eqref{eq:finlowerbd} establishes a lower bound on $\Phi_{k+1}^\star$
in terms of $F(\prm^{k+1})$ and  ${\bm e}^k$.\vspace{.2cm}

\textbf{Proving Proposition~\ref{prop:pb_aciag}}. Finally, summing up \eqref{eq:finlowerbd} and 
\eqref{eq:motivate} gives: 
\beq \begin{split}
\gap{k+1} & \leq 2 \rho^k \gap{1} + \sum_{\ell=1}^k \rho^{k-\ell} \Big( 
\sqrt{\mu \gamma} \langle {\bm e}^\ell, \prm^\star - \eprm^\ell \rangle \\
& \hspace{1.5cm} + \rho
\langle {\bm e}^\ell, \prm^\ell - \eprm^\ell \rangle 
+ \gamma \|{\bm e}^\ell \|^2 - \frac{\mu}{2} \frac{1-\mu \gamma}{\sqrt{\mu \gamma}} \| \eprm^\ell - \prm^\ell \|^2 \Big) \\
& = 2 \rho^k \gap{1} + 
\sum_{\ell=1}^k \rho^{k-\ell} \Big( 
\sqrt{\mu \gamma} \langle {\bm e}^\ell, \prm^\star - \prm^\ell \rangle \\
& \hspace{1.5cm} + 
\langle {\bm e}^\ell, \prm^\ell - \eprm^\ell \rangle 
+ \gamma \|{\bm e}^\ell \|^2 - \frac{\mu}{2} \frac{1-\mu \gamma}{\sqrt{\mu \gamma}} \| \eprm^\ell - \prm^\ell \|^2 \Big) \eqs.
\end{split}
\eeq
Let us take a look at the last summands in the above inequality: for any $ \ell \geq 1$,
\beq \begin{split}
& \sqrt{\mu \gamma} \langle {\bm e}^\ell, \prm^\star - \prm^\ell \rangle + 
\langle {\bm e}^\ell, \prm^\ell - \eprm^\ell \rangle 
+ \gamma \|{\bm e}^\ell \|^2 - \frac{\mu}{2} \frac{1-\mu \gamma}{\sqrt{\mu \gamma}} \| \eprm^\ell - \prm^\ell \|^2 \\
& \overset{(a)}{\leq} \sqrt{\mu \gamma} \| {\bm e}^\ell \| \| \prm^\star - \prm^\ell \| + \Big( 
\gamma + \frac{ \sqrt{\gamma / \mu} }{ 1 - \mu \gamma } \Big) \| {\bm e}^\ell \|^2 - \frac{\mu}{4} \frac{1-\mu \gamma}{\sqrt{\mu \gamma}} \| \eprm^\ell - \prm^\ell \|^2 \\
& \overset{(b)}{\leq} \sqrt{2 \gamma \gap{\ell}} \| {\bm e}^\ell \| + \Big( 
\gamma + \frac{ \sqrt{\gamma / \mu} }{ 1 - \mu \gamma } \Big) \| {\bm e}^\ell \|^2 - \frac{\mu}{4} \frac{1-\mu \gamma}{\sqrt{\mu \gamma}} \| \eprm^\ell - \prm^\ell \|^2 \\
& \overset{(c)}{\leq} \sqrt{2 \gamma \gap{\ell}} \| {\bm e}^\ell \| + 
\sqrt{\frac{9\gamma}{\mu}} \| {\bm e}^\ell \|^2 - \frac{\mu}{4} \frac{1-\mu \gamma}{\sqrt{\mu \gamma}} \| \eprm^\ell - \prm^\ell \|^2 \eqs,
\end{split}
\eeq
where (a) is resulted from 
the fact $\langle {\bm e}^\ell, \prm^\ell - \eprm^\ell \rangle
\leq (1/2) ( \| {\bm e}^\ell \|^2 / c + c \|  \prm^\ell - \eprm^\ell \|^2 )$ for any 
$c > 0$ 
and we have set 
$c = \frac{\mu}{2} \frac{1-\mu \gamma}{\sqrt{\mu \gamma}} $
therein; (b) is due to the relation $\| \prm^\ell - \prm^\star \| \leq \sqrt{2 \gap{\ell} / \mu}$; 
(c) is due to 
$\gamma + \frac{ \sqrt{\gamma / \mu} }{ 1 - \mu \gamma } \leq 3 \sqrt{ \gamma / \mu }$,
which can be verified through replacing $\gamma$ by its upper bound $1/(2L)$ 
in the denominator of the fraction on the left-hand-side. 
Combining the two equations above yields the desired result of Proposition.

\subsection{Proof of the Equality} \label{app:equality}
We prove 
${\bm v}^k - \eprm^k = (\sqrt{\mu \gamma})^{-1} \big( \eprm^k - \prm^k \big)$ using induction on $k$. Clearly, the said equality holds for $k=1$
since ${\bm v}^1 = \prm^1 = \eprm^1$, 
and we assume that it holds up to $k$. Consider: 
\beq \notag \begin{split}
& {\bm v}^{k+1} - \eprm^{k+1} = 
\rho {\bm v}^k + \sqrt{ \mu \gamma } \eprm^k - \sqrt{\frac{\gamma}{\mu}} {\bm g}^k - \eprm^{k+1} \\
& =\rho ( {\bm v}^k - \eprm^k ) + \eprm^k - \sqrt{\frac{\gamma}{\mu}} {\bm g}^k - \eprm^{k+1} = \frac{ \rho }{ \sqrt{\mu \gamma} } ( \eprm^k - \prm^k ) 
+ \eprm^k - \sqrt{\frac{\gamma}{\mu}} {\bm g}^k - \eprm^{k+1} \eqs,
\end{split}
\eeq
where we have used the induction hypothesis. Furthermore, 
using $\prm^{k+1} = \eprm^k - \gamma {\bm g}^k$, 
\beq
\begin{split}
& {\bm v}^{k+1} - \eprm^{k+1}  = \sqrt{\mu \gamma}^{-1} \Big( \rho (\eprm^k - \prm^k) + \sqrt{\mu \gamma} ( \eprm^k - \eprm^{k+1} )  - \gamma {\bm g}^k \Big) \\
& \overset{(a)}{=} \sqrt{ \mu \gamma }^{-1} \Big( \sqrt{\mu \gamma} ( \prm^{k+1} - \eprm^{k+1} ) + \rho (\prm^{k+1} - \prm^k ) \Big) 
= \sqrt{ \mu \gamma}^{-1} \big( \eprm^{k+1} - \prm^{k+1} \big) \eqs,
\end{split}
\eeq
where (a) is due to $\rho (\prm^{k+1} - \prm^k ) = (1 + \sqrt{\mu \gamma} ) ( \eprm^{k+1} - \prm^{k+1} )$.

\section{Proof of Proposition \ref{prop:err}} \label{app:errB}
We begin by observing that due to the $L_{H,i}$-Lipschitz continuity of the Hessian of $f_i$
and using Lemma~\ref{lem:nes},
we have:
\beq \label{eq:err_1st_eq} \begin{split}
& \| {\bm e}_{\sf ACIAG}^\ell \| = \| {\bm g}_{\sf ACIAG}^\ell - \grd F( \eprm^\ell ) \|
\leq \sum_{i=1}^m \frac{L_{H,i}}{2} \| \eprm^\ell - \eprm^{\tau_i^\ell} \|^2 \eqs.
\end{split}
\eeq
Now, expanding the right hand side of \eqref{eq:err_1st_eq} gives:
\beq  \label{eq:firstbd_e}
\begin{split}
\|{\bm e}_{\sf ACIAG}^\ell \| & \leq \sum_{i=1}^m \frac{L_{H,i}}{2} \Big\| \eprm^\ell - \eprm^{\tau_i^\ell} \Big\|^2  \leq \sum_{i=1}^m \frac{L_{H,i}}{2} \!~ \underbrace{( \ell - \tau_i^\ell )}_{\leq K} \sum_{j=\ell-\tau_i^\ell}^{\ell-1} \| \eprm^{j+1} - \eprm^j \|^2  \\
& \hspace{-1.2cm} \leq \frac{K L_{H}}{2} \hspace{-.1cm} \sum_{j=( \ell-K )_{++}}^{\ell-1} \hspace{-.1cm} \| \eprm^{j+1} - \eprm^j \|^2
= \frac{K L_{H}}{2} \hspace{-.1cm} \sum_{j=( \ell-K )_{++}}^{\ell-1} \hspace{-.1cm} \| \gamma {\bm g}_{\sf ACIAG}^j + \underbrace{\alpha ( \prm^{j+1} - \prm^j)}_{= \eprm^{j+1} - \prm^{j+1}} \|^2 \\
& \hspace{-1.2cm} \leq \frac{3 K L_H}{2} \sum_{j=( \ell-K )_{++}}^{\ell-1} \Big( \gamma^2 \big( \| {\bm e}^j \|^2 + \| \grd F( \eprm^j ) \|^2 \big) + \| \eprm^{j+1} - \prm^{j+1} \|^2 \Big) \eqs. 
\end{split}
\eeq
Remarkably,  the above bound resembles that of 
Proposition~\ref{prop:err_ciag} with the exception of the last term
that depends on $\eprm^{j+1} - \prm^{j+1}$. This is included to 
account for the extrapolated iterates used in the {\sf A-CIAG} method.

To find an upper bound of $\| {\bm e}_{\sf ACIAG}^\ell \|$ to corroborate 
Proposition \ref{prop:err}, 
in what follows, 
we will upper bound $\| {\bm e}_{\sf ACIAG}^j \|^2$ 
and $\|\grd F(\eprm^j) \|^2$, respectively. Firstly,
\beq \label{eq:ebd_1}
\begin{split}
\| {\bm e}_{\sf ACIAG}^j \| & \leq \sum_{i=1}^m \frac{L_{H,i}}{2} \Big\| \eprm^j - \eprm^{\tau_i^j} \Big\|^2 \\
& \leq \sum_{i=1}^m L_{H,i} \Big( (1+\alpha)^2 \| \prm^j - \prm^{\tau_i^j} \|^2 + \alpha^2 \| \prm^{j-1} - \prm^{\tau_i^j-1} \|^2 \Big) \eqs.
\end{split}
\eeq
Noticing that as $\| \prm^j - \prm^{\tau_i^j} \|^2 \leq 2 ( \| \prm^j - \prm^\star \|^2 + \| \prm^{\tau_i^j} - \prm^\star \|^2 ) \leq (4/\mu) ( \gap{j} + \gap{\tau_i^j} )$, 
it follows from \eqref{eq:ebd_1} that
\beq
\begin{split}
\| {\bm e}_{\sf ACIAG}^j \| & \leq \frac{4}{\mu} \sum_{i=1}^m L_{H,i} \Big( (1+\alpha)^2 ( \gap{j} + \gap{\tau_i^j} ) + \alpha^2 ( \gap{j-1} + \gap{\tau_i^j - 1} ) \Big) \\
& \leq \frac{ 8 L_H }{\mu} \big( (1+\alpha)^2 + \alpha^2 \big) \hspace{-.2cm} \max_{ (j- K-1)_{++} \leq q \leq j } \hspace{-.2cm} \gap{q} \leq \frac{ 40 L_H }{\mu} \hspace{-.2cm} \max_{ (j- K-1)_{++} \leq q \leq j } \hspace{-.2cm} \gap{q} \eqs,
\end{split}
\eeq
which implies  
\beq \label{eq:bd_E}
\begin{split}
\sum_{j=(\ell-K)_{++}}^{\ell-1} \| {\bm e}_{\sf ACIAG}^j \|^2 
& \leq 
K \Big(\frac{ 40 L_H }{\mu}\Big)^2 \max_{ (\ell - 2K-1)_{++} \leq q \leq \ell } ( \gap{q} )^2 \eqs.
\end{split}
\eeq
Secondly,
\beq \begin{split}
\| \grd F( \eprm^j ) \|^2 &
\leq 2L^2 \big( \|  \prm^j - \prm^\star \|^2 + \| \prm^j - \prm^{j-1} \|^2 \big) \leq \frac{4L^2}{\mu} \Big(  3 h^{(j)} +  2 h^{(j-1)} ) \Big) \eqs,
\end{split}
\eeq
thus
\beq \label{eq:bd_F}
\begin{split}
\sum_{j=(\ell-K)_{++}}^{\ell-1} \| \grd F( \eprm^j ) \|^2 &
\leq
\frac{20L^2 K }{\mu} \max_{ (\ell - K - 1)_{++} \leq q \leq \ell -1} h^{(q)} \eqs.
\end{split}
\eeq
Substituting \eqref{eq:bd_E} and \eqref{eq:bd_F} 
into the right hand side of \eqref{eq:firstbd_e} 
verifies Proposition~\ref{prop:err}.

\section{Step 3 in the Proof of Theorem~\ref{thm:main}} \label{app:step3_aciag}
To proceed with the proof, 
let us define the following quantity:
\beq \notag
\begin{split}
& \tilde{E}^{(\ell)} \eqdef 
\gamma^{\frac{5}{2}} \sqrt{\frac{9}{2}}  K^2 L_H \Big( 
\big( \frac{40L_H}{\mu} \big)^2 \hspace{-.3cm} \max_{ (\ell-2K-1)_{++} \leq q \leq \ell} \hspace{-.3cm} (\gap{q})^2 + \frac{20L^2}{\mu} \hspace{-.3cm} \max_{ (\ell-K-1)_{++} \leq q \leq \ell } \hspace{-.3cm} \gap{q} \Big) \\
& + 
\gamma^{\frac{9}{2}} \frac{ 81 K^4 L_H^2 }{4 \sqrt{\mu}}
\Big( 
\big( \frac{40L_H}{\mu} \big)^4 \hspace{-.3cm} \max_{ (\ell-2K-1)_{++} \leq q \leq \ell} \hspace{-.3cm} (\gap{q})^4 + \big( \frac{20L^2}{\mu} \big)^2 \hspace{-.3cm} \max_{ (\ell-K-1)_{++} \leq q \leq \ell }  \hspace{-.3cm} (\gap{q})^2 
\Big) \eqs.
\end{split}
\eeq
Using 
Proposition \ref{prop:err}, we obtain:
\beq \label{eq:thm1_inter}
\begin{split}
& \sqrt{2 \gamma \gap{\ell}} \| {\bm e}_{\sf ACIAG}^\ell \| + \sqrt{\frac{9\gamma}{\mu}} \| {\bm e}_{\sf ACIAG}^\ell \|^2  \\
&  \leq \tilde{E}^{(\ell)} + \hspace{-.4cm} \sum_{j=(\ell-K+1)_{++}}^{\ell} \hspace{-.5cm} \Big( \sqrt{\frac{9 \gamma \gap{\ell} K^2 L_H^2}{2}} \| \prm^{j} - \eprm^{j} \|^2 + \frac{27 K^3 L_H^2}{4} \sqrt{\frac{9\gamma}{\mu}} 
\| \prm^j - \eprm^j \|^4 \Big)  \eqs.
\end{split}
\eeq
We need to further bound $\gap{k}$ [recall for \eqref{eq:step1} in Proposition \ref{prop:pb_aciag}]
in terms of itself to create a `recursion' for $\gap{k}$. 
To upper bound the right hand side of \eqref{eq:step1}, let us start from 
\eqref{eq:thm1_inter}. It follows that
\beq \notag \label{eq:exactbd}
\begin{split}
& \sum_{\ell=1}^k \rho^{k-\ell} \Big( \sqrt{2 \gamma \gap{\ell}} \| {\bm e}^\ell \| + \sqrt{\frac{9\gamma}{\mu}}  \| {\bm e}^\ell \|^2 - \frac{\mu}{4} \frac{1-\mu \gamma}{\sqrt{\mu \gamma}} \| \eprm^\ell - \prm^\ell \|^2 \Big) \leq \sum_{\ell=1}^k \rho^{k-\ell}
\bigg(
\tilde{E}^{(\ell)} + \\
& \Big( 
\hspace{-.1cm} \sum_{j=\ell}^{\min\{k,\ell+K-1\}} \hspace{-.5cm} \Big( 
\sqrt{\frac{9 \gamma K^2 L_H^2  \gap{j}}{2}  } + \frac{81 K^3 L_H^2}{4} \sqrt{\frac{\gamma}{\mu}}  \| \prm^\ell - \eprm^\ell \|^2 \Big) 
- \frac{\mu}{4} \frac{1-\mu \gamma}{\sqrt{\mu \gamma}} \Big) \| \prm^\ell - \eprm^\ell \|^2
\bigg).
\end{split}
\eeq
Moreover, we observe   for $\ell \geq 2$:
\beq
\| \prm^\ell - \eprm^\ell \|^2 \leq 2 ( \| \prm^\ell - \prm^\star \|^2 + \| \prm^{\ell-1} - \prm^\star \|^2 )
\leq \frac{4}{\mu} \big( \gap{\ell} + \gap{\ell-1} \big)  \eqs,
\eeq
The coefficient in front of the last $\| \prm^\ell - \eprm^\ell \|^2$ 
term can   be upper bounded as:
\beq \notag
\tilde{C}^{(\ell,k)} \eqdef 
\gamma K^2 L_H  \sqrt{\frac{9}{2}} \hspace{-.1cm} \max_{ \ell \leq q \leq \min\{ \ell+K-1,k \}} \hspace{-.3cm} (\gap{q})^{\frac{1}{2}}
+ {\gamma} \frac{81 K^4 L_H^2}{\mu^{\frac{3}{2}}}
\big( \gap{\ell} + \gap{\ell-1} \big)  - 
\frac{\mu}{4} \frac{1-\mu \gamma}{\sqrt{\mu}}.
\eeq
If we define 
\beq \label{eq:Edef} \begin{split}
& E^{(\ell,k)} \eqdef \tilde{E}^{(\ell)}  + \tilde{C}^{(\ell,k)} \frac{\| \prm^\ell - \eprm^\ell \|^2}{\sqrt{\gamma}} \eqs,
\end{split}
\eeq
where $E^{(\ell,k)} = E^{(\ell,k -1)}$ for all $k \geq \ell + m$. 
Applying Proposition~\ref{prop:pb_aciag} readily shows 
\beq \label{eq:exact_aciag}
\gap{k+1} \leq 2 ( 1 - \sqrt{\mu \gamma} )^k \gap{1} + \sum_{\ell=1}^k (1 - \sqrt{\mu \gamma})^{k - \ell} E^{(\ell,k)} \eqs.
\eeq

\textbf{Concluding the Proof of Theorem~\ref{thm:main}}.
Our goal is to analyze \eqref{eq:exact_aciag} using Proposition~\ref{prop:gen_aciag}. Let us recognize that:
\[
\fun{k} = \bgap{k},~p = (1-\sqrt{\mu \gamma}),~b = 2,~M= 2K+1,~\eta_1 = \frac{3}{2},~\eta_2 = \frac{5}{2}, \eta_3 = 2,~\eta_4 = 4 
\]
\beq \notag
\begin{split}
& s_1 = \gamma^{\frac{5}{2}} \sqrt{\frac{9}{2}} K^2 L_H \frac{20L^2}{\mu},~
s_2 = \gamma^{\frac{5}{2}} \sqrt{\frac{9}{2}} K^2 L_H \big( \frac{40L_H}{\mu} \big)^2, \\
& s_3 = \gamma^{\frac{9}{2}} \frac{81 K^4 L_H^2}{4\sqrt{\mu}} \big(\frac{20L^2}{\mu}\big)^2,~
s_4 = \gamma^{\frac{9}{2}} \frac{81 K^4 L_H^2}{4\sqrt{\mu}} \big( \frac{40L_H}{\mu} \big)^4 \eqs,
\end{split}
\eeq
\[
c = \frac{\mu}{4} \frac{1 - \mu \gamma}{\sqrt{\mu}},~D^{(\ell)} = \frac{ \| \prm^\ell - \eprm^\ell \|^2 }{\sqrt{\gamma}},~f( \bgap{q} ) = \gamma \Big( K^2 L_H \sqrt{\frac{9}{2}} (\bgap{q})^{\frac{1}{2}} 
+ \frac{162 K^4 L_H^2}{\mu^{\frac{3}{2}}} \bgap{q} \Big) \eqs.
\]
The conditions in \eqref{eq:condition} are satisfied when
\beq
\begin{split}
& \frac{\sqrt{\mu}}{4} - \gamma \Big( K^2 L_H \sqrt{9} (\bgap{1})^{\frac{1}{2}} 
+ \frac{324 K^4 L_H^2}{\mu^{\frac{3}{2}}} \bgap{1} + \frac{\mu^{\frac{3}{2}}}{4} \Big) \geq 0 \\
& \Longleftrightarrow \gamma \leq \frac{\sqrt{\mu}}{4}
\Big( K^2 L_H \sqrt{9} (\bgap{1})^{\frac{1}{2}} 
+ \frac{324 K^4 L_H^2}{\mu^{\frac{3}{2}}} \bgap{1} + \frac{\mu^{\frac{3}{2}}}{4} \Big)^{-1}
\eqdef \frac{\bar{c}_3}{L} \eqs,
\end{split}
\eeq
and
\beq \begin{split}
1  > (1-\sqrt{\mu \gamma}) & + \gamma^{\frac{5}{2}} \sqrt{\frac{9}{2}} K^2 L_H \Big( \frac{20L^2}{\mu} (2 \bgap{1})^{\frac{1}{2}} +  \big( \frac{40L_H}{\mu} \big)^2 (2 \bgap{1})^{\frac{3}{2}} \Big) \\
& + \gamma^{\frac{9}{2}} \frac{81 K^4 L_H^2}{4\sqrt{\mu}} \Big( \big(\frac{20L^2}{\mu}\big)^2 (2 \bgap{1} ) +  \big( \frac{40L_H}{\mu} \big)^4 (2 \bgap{1})^3 \Big) \eqs,
\end{split}
\eeq
that can be implied by
\beq
\begin{split}
& \gamma < \left( \frac{\sqrt{\mu}}{\sqrt{18} K^2 L_H}\Big( \frac{20L^2}{\mu} (2 \bgap{1})^{\frac{1}{2}} +  \big( \frac{40L_H}{\mu} \big)^2 (2 \bgap{1})^{\frac{3}{2}} \Big)^{-1}  \right)^{\frac{1}{2}} \eqdef \frac{\bar{c}_1}{L}~~~~\text{and} \\
& \gamma <  \left(  \frac{2 {\mu}}{81 K^4 L_H^2} 
\Big( \big(\frac{20L^2}{\mu}\big)^2 (2 \bgap{1} ) +  \big( \frac{40L_H}{\mu} \big)^4 (2 \bgap{1})^3 \Big)^{-1} \right)^{\frac{1}{4}} \eqdef \frac{\bar{c}_2}{L} \eqs.
\end{split}
\eeq
Substituting these constants into Proposition~\ref{prop:gen_aciag} proves the claims in Theorem~\ref{thm:main}.
 
\section{Proof of Proposition \ref{prop:gen_aciag}} \label{app:prop3}
Define $\{ \bfun{k} \}_{k \geq 1}$ that satisfies:
\beq \label{eq:aux_sys}
\bfun{k+1} = p^k b \bfun{1} + 
\sum_{\ell=1}^k p^{k-\ell} \Big( \sum_{j=1}^J s_j \max_{ (\ell- M)_{++} \leq q \leq \ell } (\bfun{q})^{\eta_j} \Big),~~\bfun{1} = \fun{1} \eqs,
\eeq 
By subtracting $p \bfun{k}$ from $\bfun{k+1}$, 
\eqref{eq:aux_sys} can be alternatively expressed as:
\beq \label{eq:aux_sys_e}
\bfun{k+1} - p \bfun{k} = \sum_{j=1}^J s_j \max_{ (k- M)_{++} \leq q \leq k } (\bfun{q})^{\eta_j}  \eqs.
\eeq
Now, consider the statements (i) and (ii) in \eqref{eq:statement_aciag} 
as
the following event:
\beq \notag \label{eq:induction}
\begin{split}
{\cal E}_z = \big\{ & \!~ \bfun{(z-1)M + k+1} \geq \fun{(z-1)M + k+1},   ~\bfun{(z-1)M + k+1} \leq \delta^z (b \bfun{1} ),~
k = 1,..., M
\big\} \eqs, 
\end{split}
\eeq
for all $z \geq 1$. 
We shall prove that ${\cal E}_z$ is true for $z=1,2,...$ using induction.\vspace{.2cm}

\textbf{Base case with $z=1$}. To prove ${\cal E}_1$, let us apply another induction on $k$
inside the event. 
For the base case of $k=1$, 
\beq
\begin{split}
\bfun{2} 
& \geq  p ( b \fun{1} ) + \sum_{j=1}^J s_j (\fun{1})^{\eta_j} - ( \bar{f} - f( \fun{1})) D^{(1)} = \fun{2} \eqs,
\end{split}
\eeq
where we used the fact $\bar{f} \geq f( b \fun{1} ) \geq f( \fun{1} )$. Furthermore, 
the base case holds as:
\beq \textstyle
\bfun{2} = (b \bfun{1}) \Big( p + (1/b) \sum_{j=1}^J s_j ( \bfun{1} )^{\eta_j - 1} \Big) 
\leq \delta ( b \bfun{1} ) \eqs.
\eeq

For the induction step, 
suppose that the statements in \eqref{eq:induction} are also true up to $k=k' - 1$ with $z=1$
such that $\bfun{k'} \geq \fun{k'}$ and $\bfun{k'} \leq \delta ( b \bfun{1} )$. 
Consider the case of $k=k'$,
we observe that $\bar{f} \geq f( b \fun{1} ) \geq f (\delta b \fun{1} ) \geq f( \bfun{q} ) \geq f( \fun{q} )$ for all $q=1,...,k'$. Therefore, we can lower bound
$\bfun{k'+1}$ as:
\beq \notag \begin{split}
&\bfun{k'+1} = p^{k'} ( b \bfun{1} ) + \sum_{\ell=1}^{k'} p^{k'-\ell} \Big( \sum_{j=1}^J s_j \max_{ (\ell-M)_{++} \leq q \leq \ell} (\bfun{q})^{\eta_j} \Big) \\
& \geq  p^{k'} ( b \fun{1} ) + \sum_{\ell=1}^{k'} p^{k'-\ell} \Big( \sum_{j=1}^J s_j \max_{ (\ell-M)_{++} \leq q \leq \ell} (\fun{q})^{\eta_j}  - \big( \bar{f} - \max_{\ell \leq q \leq k'} f(\fun{q})  \big) V^{(\ell)}  \Big) , 
\end{split}
\eeq
where the right hand side is exactly $\fun{k'+1}$; 
also, using \eqref{eq:aux_sys_e}, we can show:
\beq
\begin{split}
\bfun{k'+1} 
& \leq   ( b \bfun{1} ) \Big( \delta p + \sum_{j=1}^J s_j (b \bfun{1})^{\eta_j-1} \Big) 
\leq \delta ( b \bfun{1} ) \eqs.
\end{split}
\eeq

\textbf{Induction Case}.
For the induction case, suppose that ${\cal E}_z$ is true 
for all $z$ up to $z'$. We consider the case when $z = z' + 1$.
Once again, we   apply another induction on $k$.  
In the base case of $k = 1$ and $z=z' + 1$, we have
\beq \notag \begin{split}
& \bfun{z'M+2} = p^{z'M+1} ( b \bfun{1} ) + 
\sum_{\ell=1}^{z'M+1} p^{z'M+1-\ell} \Big( \sum_{j=1}^J s_j \max_{ (\ell-M)_{++} \leq q \leq \ell} (\bfun{q})^{\eta_j} \Big) \\
& \geq p^{z'M+1} ( b \fun{1} ) + 
\sum_{\ell=1}^{z'M+1} p^{z'M+1-\ell} \Big( \sum_{j=1}^J s_j \max_{ (\ell-M)_{++} \leq q \leq \ell} (\fun{q})^{\eta_j} \\
& \hspace{5cm}
- \big( \bar{f} - \max_{\ell \leq q \leq z'M + 1} f(\fun{q})  \big) V^{(\ell)} \Big) = \fun{z'M+2} \eqs,
\end{split}
\eeq
where we used $\bar{f} \geq f( b \fun{1} ) \geq f ( \bfun{q} ) \geq f( \fun{q} )$
for all $q$ up to $q = z'M+1$ (by the induction hypothesis).  
Furthermore, the base case holds since:
\beq \begin{split}
\bfun{z'M+2} & = p \bfun{z'M+1} + \sum_{j=1}^J s_j \max_{ (z'M+1-M)_{++} \leq q \leq z'M+1 } ( \bfun{q} )^{\eta_j} \\
& \leq \delta^{z'} (b \bfun{1}) \Big( p + \sum_{j=1}^J s_j (\delta^{z'})^{\eta_j-1} (b \bfun{1})^{\eta_j-1} \Big) \leq \delta^{z'+1} ( b \bfun{1} ) \eqs.
\end{split} 
\eeq

Let the statements in ${\cal E}_z$ be true up to $k=k' - 1$, $z=z'+1$. With $k = k'$, 
\beq \notag
\begin{split}
\bfun{ z'M + k' + 1 } 
& \geq p^{z'M+k'} ( b \fun{1} ) + 
\sum_{\ell=1}^{z'M+k'} p^{z'M+k'-\ell} \Big( \sum_{j=1}^J s_j \max_{ (\ell-M)_{++} \leq q \leq \ell} (\fun{q})^{\eta_j} \\
& \hspace{3cm} - \big( \bar{f} - \max_{\ell \leq q \leq z'M + k'} f(\fun{q})  \big) V^{(\ell)}  \Big)  = \fun{z'M + k' + 1} \eqs,
\end{split}
\eeq
\beq \begin{split}
\bfun{z'M+k'+1} 
& \leq \delta^{z'} (b \bfun{1}) \Big( \delta p + \sum_{j=1}^J s_j (\delta^{z'})^{\eta_j-1} (b \bfun{1})^{\eta_j-1} \Big) \leq \delta^{z'+1} ( b \bfun{1} ) \eqs.
\end{split} 
\eeq
The induction case is thus proven.
This shows that the event ${\cal E}_z$ is true for all $z \geq 1$.\vspace{.2cm}

\textbf{Proving Statement (iii)}. We apply statement (ii) to prove (iii). 
From \eqref{eq:aux_sys_e}, 
\beq
\begin{split}
\frac{ \bfun{k+1} }{ \bfun{k} } & = p + \frac{1}{ \bfun{k} } \sum_{j=1}^J s_j \max_{ (k-M)_{++} 
\leq q \leq k} (\bfun{q} )^{\eta_j} \eqs.
\end{split}
\eeq
For any $q \in [(k-M)_{++}, k]$, we have 
\beq
\frac{ (\bfun{q})^{\eta_j} }{\bfun{k}} = \frac{ \bfun{q} }{ \bfun{k} } (\bfun{q})^{\eta_j - 1}
\leq \frac{ \bfun{q} }{ \bfun{k} } \Big( \delta^{\lceil (q-1) / M \rceil} ( b \fun{1} ) \Big)^{\eta_j - 1} \eqs.
\eeq
Since $\eta_j > 1$ and $|q-k| \leq M$, we have $\delta^{\lceil (q-1) / M \rceil ( \eta_j - 1 )}  \rightarrow 0$
as $k \rightarrow \infty$, 
moreover as $\bfun{k+1} / \bfun{k} \geq p$ for all $k \geq 1$,  
$\bfun{q} / \bfun{k} \leq p^{-M}$ for all $q$. 
Therefore, we get 
\beq
\lim_{ k \rightarrow \infty } \frac{ \max_{ (k-M)_{++} \leq q \leq k} (\bfun{q} )^{\eta_j}  }{ \bfun{k} } 
= 0,~\forall~j 
\Longrightarrow \lim_{ k \rightarrow \infty } \frac{ \bfun{k+1} }{ \bfun{k} } = p \eqs.
\eeq

%
%

\end{document}